\documentclass[a4paper] {article}

\usepackage{amsmath,amssymb, amsthm,eucal}
\usepackage{amscd}
\usepackage{t1enc}
\usepackage[all]{xy}

\theoremstyle{plain} 
\newtheorem{proposition}{Proposition}[section]  
\newtheorem{lemma}[proposition]{Lemma} 
\newtheorem{corollary}[proposition]{Corollary} 
\newtheorem{theorem}[proposition]{Theorem} 

\theoremstyle{definition} 
\newtheorem{definition}[proposition]{Definition} 
\newtheorem{fact}[proposition]{Fact}
\newtheorem{remark}[proposition]{Remark} 
\newtheorem{terminology}[proposition]{Terminology}
\newtheorem{example}[proposition]{Example} 
\newtheorem{recall}[proposition]{Recall}
\newtheorem{prel}[proposition]{Preliminaries and notations}
\newtheorem{question}[proposition]{Question}
\newtheorem{observation}[proposition]{Observation}

\newtheorem{questions}[proposition]{Questions}
\newtheorem{observations}[proposition]{Observations}

\newtheorem{remarks}[proposition]{Remarks}
\newtheorem{notation}[proposition]{Notation}

\def\wl{\vspace{\baselineskip}}

\makeatletter
\renewcommand{\section}{\@startsection{section}{1}{0cm}{-1.5\baselineskip}%
{0.5\baselineskip}{\centering\large\rmfamily\upshape\bfseries}}
\renewcommand{\subsection}{\@startsection{subsection}{2}{0cm}{-1\baselineskip}%
                   {0.2\baselineskip}{\bfseries}}
\def\@seccntformat#1{\csname the#1\endcsname.~~}
\makeatother

\begin {document}

\title{Complete intersections of dimension zero : variations on a theme of Wiebe}
\author{Anne-Marie Simon and Jan R. Strooker}

\date{September 7 2006}

\maketitle

\begin{abstract}Wiebe's criterion, which recognizes complete intersections of dimension zero among the class of noetherian local rings, is revisited and exploited in order to provide information on what we call C.I.0-ideals (those such that the corresponding quotient is a complete intersection of dimension zero) and also on chains of C.I.0-ideals.

A correspondence is established between C.I.0-ideals and a certain kind of matrices which we call $x$-nice, and a chain of C.I.0-ideals corresponds to a factorization of some $x$-nice matrix.

When the local ring $A$ itself is a complete intersection of dimension zero, a C.I.0-ideal is necessarily of the form $(0:bA)$ for some $b\in A$. Some criteria are provided to recognize whether an ideal $(0:bA)$ is C.I.0 or not. When $y$ is a minimal generator of the maximal ideal of $A$, it is also proved that the ideals $yA$ and $(0:yA)$ are C.I.0 simultaneously and  that this is the case exactly when the ideal $(0:yA)$ is principal. These C.I.0-ideals of the form $(0:yA)$, $y$ being a minimal generator of the maximal ideal, are investigated. They are of interest because the smallest nonnull C.I.0-ideal in a  strict chain of C.I.0-ideals of the maximal length is necessarily of that form, and their existence has some implications for a realization of the ring, i.e. for the way the ring can can be viewed as a quotient of a regular local ring.

\wl MSC: 13E10, 13H10.

\wl Keywords: Complete intersection, Gorenstein, Artinian local ring, Annihilator, Chain.
\end{abstract}

\section*{Introduction}
A noetherian local ring is called a complete intersection if its completion with respect to its maximal-adic topology is the quotient of a regular local ring by an ideal generated by a regular sequence. Complete intersections are Gorenstein and complete intersections of dimension zero are recognised by Wiebe's criterion.
 
Here we are interested in the proper ideals of a noetherian local ring $A$ such that the corresponding quotient is a complete intersection of dimension zero and call them zero-dimensional complete intersections ideals (C.I.0-ideals). How can we produce them, how can we produce the elements of $A$, the annihilator of which is a C.I.0-ideal, and what can we say about chains of C.I.0-ideals ? It turns out that Wiebe's criterion, together with the arguments of its proof, also produces some answers to these questions.

This criterion states that a noetherian local ring $A$ with a sequence $x=(x_1,\cdots ,x_n)$ generating its maximal ideal is a complete intersection of dimension zero if and only if there is a matrix $\psi \in A^{n\times n}$ such that the row-matrix $x\cdot \psi$ is null while $\text{det}(\psi)$ is nonnull. When this is the case, $\text{det}(\psi)$ generates the socle of $A$. We call such a matrix $\psi$ an $x$-Wiebe matrix for $A$. With Wiebe's criterion we can see that an $x$-Wiebe matrix for a ring encodes the structure of the ring. And it turns out that the factorizations of the $x$-Wiebe matrices of $A$ correspond to the chains of C.I.0-ideals in $A$, just like the factorizations of an element generating the socle of a Gorenstein local ring of dimension zero correspond to the chains of Gorenstein ideals in that ring.

The first section introduces notations, describes linkage of completely sequent sequences in an arbitrary commutative ring with unit and contains some more facts on regular sequences in a Gorenstein local ring.

In the second section, we present a slight extension of Wiebe's criterion in the form we need. We include a  somewhat simplified proof which also allows us to keep, for later use, the information given by the arguments. Doing so we show that these C.I.0-ideals are stemming from what we call $x$-nice matrices. By an $x$-nice matrix of a noetherian local ring $A$ with a sequence $x=(x_1,\cdots,x_n)$ generating its maximal ideal, we mean a matrix $\varphi\in A^{n\times n}$ such that its determinant does not belong to the ideal generated by the entries of the row-matrix $x\cdot \varphi$. On the way, and in analogy with Wiebe's criterion, we also provide criteria to recognize when the annihilator of an element is a C.I.0-ideal. We provide many examples and small applications.

In the intermediary third section we collect some facts about chains of Gorenstein ideals in a zero-dimensional Gorenstein local ring. They correspond to factorizations of an element generating the socle $0:M$ of the ring so that  their lengths are bounded by $\text{max}\{i \mid M^i\neq 0\}$.

The fourth section is concerned with chains of C.I.0-ideals. They are treated in analogy with the Gorenstein case.

 To an inclusion $I_0\subset I_1$ of C.I.0-ideals in any ring corresponds a factorization of an $x$-nice matrix belonging to $I_0$, the first factor corresponding to $I_1$ while the determinant of the second factor generates, modulo $I_0$, the annihilator of $I_1/I_0$. 

Consequently, when the ring $A$ itself is a complete intersection of dimension zero, we can produce every C.I.0-ideal and the principal ideal of which this is the annihilator by factoring the $x$-nice matrices belonging to the null ideal (i.e. the $x$-Wiebe matrices). Moreover, and in analogy with the Gorenstein case, we have a correspondence between chains of C.I.0-ideals and factorizations of $x$-Wiebe matrices.

At the end of the section we raise some problems concerning the length of the chains of C.I.0-ideals in a zero-dimensional complete intersection $A$. We know that these lengths are bounded
 and we have an example of such an $A$ with strict saturated chains of C.I.0-ideals of different length. But when do we have a strict chain of C.I.0-ideals of the maximal length?
We show that the existence of  such a chain is related to questions concerning the minimal generators of the maximal ideal: when is their annihilator not only Gorenstein but also C.I.0 ? 
 
 The fifth section contains partial answers to the questions raised before. In a zero-dimensional complete intersection $A$ we first investigate those minimal generators $y$ of the maximal ideal, the annihilator of which is  C.I.0. It turns out that, if $y$ is a minimal generator of the maximal ideal $M$ of $A$, then the ideals $yA$ and $0:yA$ are C.I.0-ideals simultaneously, and that it is the case exactly when the ideal $0:yA$ is also principal. Moreover, the existence of such a $y$  has some implications for the realizations of $A$.

In the sixth section we provide a tentative classification of C.I.0-ideals and a description of the set of $x$-nice matrices belonging to a given C.I.0-ideal.

The last section contains another sufficient condition under which the ideal generated by a minimal generator of the maximal ideal of a complete intersection of dimension zero is itself a C.I.0-ideal.

\wl The second author would like to thank Larry Smith, Göttingen. During a workshop there, we became aware that, from different perspectives, we had come across similar problems. Larry pointed out \cite{Wi} and other references, and made a preliminary version of \cite{ms} available. In fact, proposition \ref{p1.1} is our attempt to understand proposition VI.3.1 in this Tract.

\wl We put this text in the arXiv, while not submitting it for publication now. It does not have the tight format of a research paper. It reads more like a relaxed ramble through a pleasant country side, through familiar copses and less known spinneys. We stop by the wayside to ask questions, examine examples of the vegetation and are not in much of a hurry. Also, though we do provide some insight, certain rather natural problems remain unsolved.

\section {Linkage of completely sequent sequences and applications}
\begin{prel}\label{prel} Let $a=(a_1,\cdots,a_n)$ be a sequence of elements of a commutative ring $A$ with unit. The ideal generated by this sequence will be denoted by $J(a)$ (or simply by $(a_1,\cdots ,a_n)$ in the examples where such a  matrix notation is harmless). We alternatively view $a$ as a sequence or as a row matrix. Sometimes we also identify a map $A^m \rightarrow A^n$ with its matrix $\gamma$ in the canonical bases, doing so, we view $A^n$ as a set of column matrices ($A^n \equiv A^{n\times 1}$), so $\gamma \in A^{n\times m}$.

Now let $\varphi \in A^{n\times n}$ be a square matrix of size n, so that $J(a \cdot \varphi)\subseteq J(a)$. We denote by $\text{det}(\varphi)$ the determinant of $\varphi$, by $\varphi^t$ the transpose of $\varphi$ and by $\varphi^c$ the transpose of the cofactor matrix of $\varphi$. Thus $\varphi^c$ is the adjoint matrix of $\varphi$ with $\varphi \cdot \varphi^c = \varphi^c \cdot \varphi =\text{diag}(\text{det}(\varphi))$.
From $a \cdot \varphi \cdot \varphi^c = \text{det}(\varphi) \cdot a$, we deduce that $\text{det}(\varphi) \in J(a \cdot \varphi) : J(a)$.

Much more can be said when the sequence $a \cdot \varphi = (b_1 ,\cdots, b_n) = b$ is completely sequent, which means that the Koszul homology $H_i(b, A) = 0$ for $i\geq1$, equivalently that the Koszul cohomology $H^i(b,A) = 0$ for $i < n$. In this case, we note that the sequence $a$ is also completely sequent. Indeed, for any $A$-module $W$, let $h^-(a,W) = \text{inf}\{i\mid H^i(a,W)\neq 0\}$, $h^-(a,W)$ is a natural number or $\infty$; as $J(b) \subseteq J(a)$, we have $n\leq h^-(b,W) \leq h^-(a,W)$, see \cite{si88} or use ( \cite{st90}, 5.3.11 and 6.1.6).
\end{prel}

The following is probably well-known, part of it can be found in (\cite{Wi}, Satz 2) with a somewhat different argument, in a slightly more restricted situation it can also be found in (\cite{ms}, VI.3.1).

\begin{proposition} \label {p1.1}
Let $a =(a_1,\cdots, a_n)$ and $b = (b_1, \cdots, b_n)$ be two sequences of elements of a commutative ring $A$ (with unit), such that $J(b)\subset J(a)$ and such that the sequence $b$ is completely sequent. Let $\varphi \in A^{n\times n}$ be any matrix with $b = a\cdot \varphi$.

Then the sequence $a$ is also completely sequent and
\begin{equation}
J(b) : J(a) = J(b) + \text{det}(\varphi)A ,
\end{equation}

\begin{equation}
J(b) : (J(b) + \text{det}(\varphi)A) = J(a).
\end{equation}

Moreover, if $J(a) \neq A$, then $\text{det}(\varphi \notin J(b)$.
\end{proposition}

\begin{proof}
We already know that $\text{det}(\varphi) \in J(b):J(a)$, we also know that both Koszul complexes $K.(b,A)$ and $K.(a,A)$ are exact in positive degrees (see \ref{prel}) and provide finite free resolutions of $A/J(b)$ and $A/J(a)$. Let $z \in J(b):J(a)$. Multiplication by $z$ gives a commutative diagram, where $p$ is the natural surjection
\begin{displaymath}
\begin{array}{ccc}
$A$ & \longrightarrow & $A/J(b)$\\
\parallel & {} & {\downarrow p}\\
$A$ & \longrightarrow & $A/J(a)$\\
{ z\downarrow \phantom{i}} & {} & {\downarrow f}\\
$A$ & \longrightarrow & $A/J(b)$\\
\end{array}
\end{displaymath}
We lift $p$ and $f$ into maps of Koszul complexes, $\wedge \varphi$ being a lifting of $p$ :
\begin{displaymath}
\begin{array}{cccccccccccccc}
K.(b,A): & 0 &\rightarrow  &  A &  \stackrel{b^t}{\rightarrow} & A^n & \rightarrow & \cdots & \rightarrow & A^n & \stackrel{b}{\rightarrow} & A &\rightarrow & 0\\
 \phantom{ii}\downarrow\wedge\varphi &{}  & {} &\phantom{iiiiii} \downarrow \text{det}(\varphi) & {} & \downarrow &{}&{}&{} & \phantom{}\downarrow \varphi & {} & \parallel  &{} &{} \\ 
 K.(a,A): & 0 &\rightarrow & A & \stackrel{a^t}{\rightarrow}  & A^n & \rightarrow & \cdots & \rightarrow & A^n & \stackrel{a}{\rightarrow } & A & \rightarrow & 0\\
\phantom{i}\downarrow\zeta &{}  & {} & \phantom{i}\downarrow c
 & {} & \downarrow & {} & {} & {} & \downarrow & {}  & \phantom{i}\downarrow z & {} &{} \\
K.(b,A) : & 0 &\rightarrow & A & \stackrel{b^t}{\rightarrow} & A^n & \rightarrow & \cdots &\rightarrow & A^n &\stackrel{b}{\rightarrow} & A & \rightarrow & 0\\
\end{array}
\end{displaymath}

Both maps $\zeta \circ\wedge\varphi$ and multiplication by $z$ are liftings of $f\circ p$,  they are thus homotopic and we have a map $g:A^n\longrightarrow A$ such that 
  $(z-c\circ \text{det}(\varphi)) = g\circ b^t$.
 Consequently $(z-c\cdot \text{det}(\varphi)) \in J(b)$ and $J(b):J(a) \subseteq J(b) + \text{det}(\varphi)A$. The other inclusion being already known, we have the first equality.

\wl On the other hand, we have
\begin{displaymath}
(J(b):J(a))/J(b) \simeq Hom_A(A/J(a), A/J(b)) \simeq H^0(a, A/J(b)).
\end{displaymath}
Now we split off the resolution of $A/J(b)$ given by the Koszul complex$K.(b,A)$, this give short exact sequences
$$ 0 \rightarrow A  \stackrel{(b_1,\cdots ,b_n)^t}\rightarrow A^n \rightarrow M_{n-1} \rightarrow 0, \quad \cdots \quad ,  
0 \rightarrow M_1  \rightarrow A \rightarrow A/J(b) \rightarrow 0 ,$$ 
where $M_i$ denotes the $i^{\text th}$ sysygy of $A/J(b)$.
We apply the long exact sequence of the $H^i(a, \cdot)$ to these short exact sequences. As $h^-(a, A) \geq n$ and $H^n(a, A) = A/J(a)$, we obtain isomorphisms
$$H^0(a, A/J(b)) \simeq H^1(a, M_1) \simeq \cdots \simeq H^{n-1}(a, M_{n-1})$$
and a commutative diagram with exact rows

\[ \begin{CD}
0 @>>> H^{n-1}(a, M_{n-1}) @>>> H^n(a, A) @>{H^n(a, (b_1,\cdots ,b_n)^t)}>> H^n(a, A^n) @>>> \cdots \\
@.  @. @| @| @.\\
 @.  @. A/J(a) @>{(b_1,\cdots ,b_n)^t}>> A/J(a)^n @. 
\end{CD} \]
 In this diagram, the map $H^n(a, (b_1,\cdots ,b_n)^t)$ is null since $J(b) \subseteq J(a)$, so that we have
 $$(J(b) : J(a))/J(b) \simeq H^0(a, A/J(b)) \simeq H^{n-1}(a, M_{n-1}) \simeq H^n(a,A) = A/J(a) .$$

With identity (i) in the proposition, we obtain 
\begin{displaymath}
(J(b) + \text{det}(\varphi)A)/J(b) \simeq A/J(a),
\end{displaymath}
which gives
\begin{displaymath}
J(a) =  \text{Ann}_A((J(b) + \text{det}(\varphi)A)/J(b)) =  J(b):(J(b) + \text{det}(\varphi)A) .
\end{displaymath}
This last equality implies $\text{det}(\varphi) \notin J(b)$ when $J(a) \neq A$. 
\end{proof}

This proposition can be viewed as a result in linkage theory, as it describes a particular case of linkage of perfect ideals.

\begin{remark}\label{rfit} In the situation of the above proposition, $\text{det}(\varphi)$ generates $\text{Ann}(J(a)/J(b))$ modulo $J(b)$, moreover, $\text{det}(\varphi)$ belongs to the {\bf initial Fitting ideal} of $J(a)/J(b)$ viewed as an $A$-module. This follows from the description of the initial Fitting ideal which we now recall.

In general, if $W$ is a finitely presented $A$-module and if $A^m \stackrel{f}{\rightarrow} A^n\stackrel{p}{\rightarrow} W\rightarrow 0$ is a presentation of W (i.e. an exact sequence), then the initial Fitting ideal   $\delta^A_0(W)$ of $W$ is $I_n(\theta)$, the ideal generated by the $n\times n$-minors of a matrix $\theta$ representing $f$. Equivalently, $\delta^A_0(W)$ is generated by the elements of the shape $\text{det}(g)$, where $g:A^n \longrightarrow A^n$ is a map with $p \circ g=0$.

This initial Fitting ideal of $W$ is an invariant of $W$, and we have 
$\text{Ann}(W)^n \subseteq \delta^A_0(W) \subseteq \text{Ann}(W)$, see \cite{no} . Moreover, when $W$ is an $\bar{A}$-module for some homomorphic image $\bar{A}$ of $A$, the initial Fitting ideal of $W$, viewed as an $\bar{A}$-module, is the image of $\delta^A_0(W)$ in $\bar{A}$.
\end{remark}

With these remarks we obtain the following.

\begin{corollary}\label{cf} 
In the situation of the proposition, let us denote by $\overline{(\cdot )}$ the images modulo $J(b)$. We have: $\delta_0^{\bar A}(J(a)/J(b)) = det(\overline{\varphi })\bar A$.
\end{corollary}

\begin{proof}We just saw that our initial Fitting ideal $\delta_0^{\bar A}(J(a)/J(b))$ is generated by the elements $\text{det}(\overline \varphi _1)$, where $\varphi _1 \in   A^{n \times n}$ runs among the matrices such that $\bar{a}\cdot \overline{\varphi _1} = 0$, equivalently such that $J(a\cdot \varphi _1) \subseteq J(b)$. For such a matrix $\varphi _1$, we saw in \ref{prel} that $\text{det}(\varphi _1) \in J(b):J(a)$ and the latter is $J(b)+\text{det}(\varphi)A$ by the proposition. The conclusion follows.
\end{proof}

\begin{corollary}\label{c1.2}
Let $(S,M',K)$ be a regular local ring with its maximal ideal and residue field, and let $x'=(x'_1,\cdots, x'_n)$ be a sequence generating $M'$ minimally.\newline
Let $a'=(a'_1,\cdots,a'_n)$ be a maximal regular sequence in $S$ and let $\varphi'$ be any matrix  with $a' = x'\cdot \varphi'$.
Let $A = S/J(a')$, so that $A$ is a complete intersection of dimension zero, let $M = M'/J(a')$ and let $\varphi$ be the image of $\varphi'$ in $A$. Then 
\begin{enumerate}
\item[(i)]  $socle(A) = 0 :_A M = det(\varphi)A \simeq K$,
\item[(ii)]  $det(\varphi') \notin J(x' \cdot \varphi')$,
\item[(iii)]  $\delta^A_0(M) = det(\varphi)A \neq 0$,
\item[(iv)] $det(\varphi') \in I'$ for every ideal $I'$ of $S'$ containing $J(x' \cdot \varphi')$ properly.
\end{enumerate}
\end{corollary}

\begin{terminology} Note that an artinian local ring has two natural invariants, its {\bf embedding dimension}, namely the minimal number  of generators of its maximal ideal $M$, and its {\bf exponent}, which is $\text{min} \{r \in \mathbb N_0 \mid M^r=0\}$.
\end{terminology}

\begin{corollary}\label{c1.3}
Let $A$ be a complete intersection of dimension 0 and embedding dimension n. Let $M$ be the maximal ideal of $A$. Then $M^n \neq 0$, i.e. the exponent of $A$ is strictly greater than its embedding dimension.
\end{corollary}
\begin{proof}
 We can write $A = S/J(a')$, where $S$ and $a'$ are as in corollary \ref{c1.2} of which we preserve the notations, thus $\text{dim}(S)=n$. Since the embedding dimension of $A$ is the dimension $n$ of $S$, $J(a') \subset M'^2$ and the matrix $\varphi'$ has its entries in $M'$, thus $\text{det}(\varphi) \in M^n$, but $\text{det}(\varphi) \neq 0$.
\end{proof}

Let us go back to proposition \ref{p1.1}, searching for some partial converse.
More precisely, in the situation of this proposition \ref{p1.1}, we retain that, if $J(a)\neq A$ and if the sequence $b = a\cdot \varphi$ is regular, then so is the sequence $a$ and $\text{det}(\varphi) \notin J(a\cdot \varphi)$, and we wonder  about a converse to this?

We shall see that it holds under some rather restrictive though very interesting hypothesis, see the important proposition \ref{p1.4}.

 We note that it fails in general, see the following.
 
 \begin{example}\label{exfirst} Let $A=K[[X,Y,]]$, where $K$ is a field. Put $a=(X^2,Y)$, $\varphi = \left( \begin{array}{cc}
 1&0\\
 0&X
 \end{array}\right)$. The sequence  $a$ is regular, the sequence  $a\cdot \varphi = (X^2, XY)$ is not regular and  $\text{det}(\varphi) =X \notin J(a\cdot \varphi)$. However, the primary decomposition of the ideal $J(a\cdot \varphi)$ is $(X^2, XY)= (X^2, Y) \cap (X)$ and $det(\varphi) = X $ belongs to the  primary components of minimal height  of the ideal $J(a\cdot \varphi)$.
\end{example}
 
\wl  Now we present a partial converse to proposition \ref{p1.1}. First we recall some elementary facts about linkage in a Gorenstein local ring.
 
 \begin{definition}\label{unmix} An ideal $I$ of a noetherian local ring is called {\bf purely unmixed} or  {\bf of pure height} if its primary components have all the same height.
 
 In general, the {\bf  purely unmixed ideal} $I_u$ associated to $I$ is the larger ideal which is the intersection of all the primary components of $I$ of minimal height (these are uniquely determined).
\end{definition}

The following is well-known.

\begin{lemma}\label{gor}
Let $I$ be a grade $g$ ideal in a Gorenstein local ring $A$ and let $y=(y_1, \cdots ,y_g)$ be a maximal regular sequence contained in $I$. 

Then $J(y):I$ is a purely unmixed ideal of grade $g$ and $J(y): (J(y):I)=I_u$.
\end{lemma}

\begin {proof} Recall that grade and height coincide for an ideal in a Gorenstein local ring.
Recall also that the ideal $J(y):I$ is always purely unmixed and that $I= J(y):(J(y):I)$ in the case $I$ itself is purely unmixed, these facts are well known and can be found in many papers, see for example \cite{ps} or \cite{scli}. 

If the ideal $I$ is not purely unmixed, let $I= Q_1 \cap \cdots \cap Q_k \cap Q_{k+1} \cap \cdots \cap Q_{\ell}$ be a primary decomposition of $I$, where $\text{grade}(Q_i) = g$ for $1\leq i \leq k$ and $\text{grade}(Q_j) > g$ for $k< j \leq \ell$. With these notations we have $I_u = Q_1 \cap \cdots \cap Q_k$.
 As $\text{grade}(Q_{k+1} \cap \cdots \cap Q_{\ell}) > g$, the regular sequence $y=(y_1, \cdots , y_g)$ can be extended to a regular sequence $y' = (y_1, \cdots , y_g, z)$ contained in $Q_{k+1} \cap \cdots \cap Q_{\ell}$. We then have $z\cdot I_u \subset I$, this implies  $\forall t \in J(y):I, \quad t\cdot z \cdot I_u \subset J(y)$ and also $t\cdot I_u \subset J(y)$ since the sequence $y'$ is regular. Thus $J(y):I \subset J(y):I_u$; the other inclusion being obvious we have the wanted equality.
\end{proof}

\begin{proposition}\label{gra} Let $A$ be a Gorenstein local ring and let $a=(a_1, \cdots ,a_n)$ and $b=(b_1, \cdots b_n)$ be two sequences in $A$ such that $J(b)\subseteq J(a)$ and $\text{grade}(J(b)) < \text{grade}(J(a))$  . Let $\varphi \in A^{n\times n}$ be any matrix such that $J(a\cdot \varphi) \subseteq J(b)$.

  Then $det(\varphi) \in J(b)_u$.
\end{proposition}

\begin{proof}Let $y=(y_1, \cdots ,y_g)$ be a maximal regular sequence in $J(b)$, $g = \text{grade}(J(b))$. Since $\text{det}(\varphi) \cdot J(a) \subset J(b)$ by \ref{prel}, we have 
\[ \forall t \in J(y):J(b), \quad t\cdot \text{det}(\varphi)\cdot J(a) \subseteq J(y).\]
Since $J(a)$ contains an element $z$ regular on $A/J(y)$ because $\text{grade}(J(b)) < \text{grade}(J(a))$,  we also have 
$$ \forall t \in J(y):J(b), \quad t\cdot \text{det}(\varphi) \in J(y).$$ Thus $\text{det}(\varphi) \in J(y):(J(y):J(b))$ and this last ideal is the purely unmixed ideal $J(b)_u$, see \ref{gor}.
\end{proof}

\begin{remark}\label{cfit} In the situation of proposition \ref{gra}, we also have $\delta^0_A(J(a)/J(b)) \subseteq J(b)_u$ (see the description of the initial Fitting ideal given in \ref{rfit}).
\end{remark}

\begin{corollary}\label{rec} Let $A$ be a Gorenstein local ring and let $a=(a_1, \cdots ,a_n)$ be a regular sequence in $A$. Let $b=(b_1, \cdots b_n)$ be another sequence in $A$ such that $J(b)\subseteq J(a)$ and let $\varphi \in A^{n\times n}$ be any matrix such that $b=a\cdot \varphi$.

Then the sequence $b$ is regular if and only if $det(\varphi) \notin J(b)_u$.
\end{corollary}

\begin{proof}If the sequence $b$ is regular, then the ideal $J(b)$ is purely unmixed, equal to $J(b)_u$, and $\text{det}(\varphi) \notin J(b)_u$ by proposition \ref{p1.1}.

If the sequence $b$ is not regular, then $\text{grade}(J(b)) < \text{grade}(J(a))$ and the conclusion follows from proposition \ref{gra}.
\end{proof}

\section{Wiebe's criterion}
Here is another converse to proposition \ref{p1.1}.
It is also a key argument in Wiebe's criterion (\cite{Wi}, Satz 1 and its proof), see also (\cite{BH}, theorem 2.3.16 and its proof).

\begin{proposition}\label{p1.4}
(\emph{Wiebe})Let $S$ be a regular local ring with a sequence $x = (x_1,\cdots, x_n)$ in $S$ generating its maximal ideal $M$ minimally. Let $\varphi \in S^{n\times n}$ and $x \cdot \varphi =a =(a_1,\cdots,a_n)$.

Then $det(\varphi) \notin J(x \cdot \varphi)$ if and only if the sequence $a = x \cdot \varphi$ is a maximal regular sequence in $S$.
\end{proposition}
\begin{proof}
The if part is in proposition \ref{p1.1}.

So assume that the sequence $a$ is not regular and let us then prove that $\text{det}(\varphi) \in J(x\cdot \varphi)$. \\
For every natural number $t$, we first construct inductively a sequence $a'' = (a''_1,\cdots,a''_n)$  in  $M^{t+1}$ such that the modified sequence $a' = a+a''$ is regular. Assume we already constructed $a''_1,\cdots, a''_i,  0\leq i<n$. Let $p_1,\cdots,p_s$ be the distinct associated prime ideals of $J(a'_1,\cdots,a'_i)$, enumerated in such a way that $a_{i+1} \in p_1 \cap \cdots \cap p_r$ and $a_{i+1} \notin p_{r+1} \cup \cdots \cup p_s$, $0\leq r\leq s$. These prime ideals have all height $i < n =\text{dim}S $. With prime avoidance, we pick $b \in  M^{t+1} \setminus  (p_1 \cup \cdots \cup p_s)$ and $c \in (p_{r+1} \cap \cdots \cap p_s) \setminus (p_1 \cup \cdots \cup p_r)$. We then put $a''_{i+1} = bc$ and we observe that $a'_{i+1} = a_{i+1}+bc \notin (p_1 \cup \cdots \cup p_s)$. Thus the sequence $(a'_1,\cdots,a'_i,a'_{i+1})$ is regular. We now have our sequence $a''$.

 As the sequence $a''$ has its entries in $M^{t+1}$, we have a matrix $\varphi'' \in S^{n\times n}$ with entries in $M^t$ such that $a'' = x\cdot \varphi''$. We put $\varphi' = \varphi+\varphi''$, so that $a' =x \cdot \varphi'$. As  $\varphi$ and $\varphi'$ are equivalent modulo $M^t$, so are their determinants: $\text{det}(\varphi) - \text{det}(\varphi') \in M^t$.\\
  Put $I=J(x\cdot \varphi)$, $I'=J(x\cdot \varphi')$. We have: $I'\subseteq I+M^{t+1} \subsetneq I+M^t$, the second inclusion is strict because, by assumption, $\text{dim}(S/I)>0$. Because the sequence $a' = x\cdot \varphi'$ is regular, we then have: $\text{det}(\varphi') \in I+M^t$ (\ref{c1.2} (iv)).Then we also have $\text{det}(\varphi) \in I+M^t$. As this holds for every natural number $t$, we obtain with Krull's intersection theorem what we want: $\text{det}(\varphi) \in I=J(x\cdot \varphi)$.
\end{proof}

The condition on $\varphi$ in the above proposition will play a central role in the study of C.I.0-ideals. Let us retain it.

\begin{definition}\label{d1.1}
Let A be a noetherian local ring with a sequence $x = (x_1,\cdots,x_n)$ generating its maximal ideal (not necessarily minimally). 

 A matrix $\varphi \in A^{n\times n}$
is called {\bf $x$-nice } when $\text{det}(\varphi) \notin J(x\cdot \varphi)$. An $x$-nice matrix $\varphi$ is said to {\bf belong} to the ideal $J(x\cdot \varphi)$ and two $x$-nice matrices are called {\bf equivalent} if they belong to the same ideal.
\end{definition}

We shall also need the following.

\begin{definition}\label{d1.2} 
When the local ring $A$ is a homomorphic image of a regular local ring $S$, $A \simeq S/Q$, we say that $(S,Q)$ is a realization of $A$. Let $x = (x_1,\cdots,x_n)$ be a sequence in $A$ generating its maximal ideal $M$ (not necessarily minimally).By an {\bf $x$-realization $(S,Q,x')$ of $A$ }, we mean a realization $(S,Q)$ of $A$ with a sequence $x' = (x'_1,\cdots,x'_n)$ in $S$ generating the maximal ideal $M'$ of $S$ minimally, and such that the image of $x'_i$ in $A$ is $x_i$.
\end{definition}

With  Cohen's structure theorem, every complete local ring with a sequence $x$ as above has an $x$-realization.

We are ready to present a slight extension and a more precise version of Wiebe's criterion, in the form we want for later use. In particular, we state some facts for rings of positive dimension.

\begin{theorem}\label{twi}
Let $A$ be a noetherian local ring with a sequence $x = (x_1,\cdots,x_n)$ generating its maximal ideal $M$ (not necessarily minimally).Let $I$ be an ideal in $A$, put $\bar A = A/I$ and let us denote by $\overline{( \cdot )}$ the images modulo $I$.
Then the following conditions are equivalent:
\begin{enumerate}
\item[(i)] the ring $\bar A$ is a complete intersection of dimension 0,
\item[(ii)] $\delta^{\bar A}_0 (\overline M) \neq 0$,
\item[(iii)] there is a matrix $\varphi \in A^{n\times n}$ with $J(x\cdot \varphi) \subseteq I$ and $det(\varphi) \notin I$.
\end{enumerate}

When these conditions are satisfied, let $\varphi$ be a matrix as in (iii), let $(S,Q,x')$ be an $\bar x$-realization of $\bar A$, and let $\varphi' \in S^{n\times n}$ be a matrix whose image in $\bar A$ is $\overline \varphi$. Then:
\begin{enumerate}
\item[(a)] the sequence $x'\cdot\varphi'$ is a maximal regular sequence in $S$ and $Q = J(x'\cdot \varphi')$,
\item[(b)] $0 :_{\bar A} \overline M = det(\overline\varphi)\bar A = \delta _0^{\bar A}(\overline M)$,
\item[(c)] $\varphi$ is an $x$-nice matrix belonging to $I$: $I = J(x\cdot\varphi)$ and $J(x\cdot \varphi) : M = J(x\cdot \varphi) + det(\varphi)A$.
\end{enumerate}
\end{theorem}

\begin{proof}
$(i)\Rightarrow(ii).$ \quad This is a consequence of corollary\ref{c1.2} (iii).

 $(ii)\Leftrightarrow(iii).$ \quad This follows from the description of the initial Fitting ideal given in \ref{rfit} and lifting from $\bar A$ to $A$.

 $(iii)\Rightarrow(i),(a),(b).$ \quad Let $\hat {\bar A}$  be the $\overline M$-adic completion of $\bar A$, $\bar A \subseteq \hat {\bar A}$ and let $(S,Q,x')$ be an $\bar x$-realization of  $\hat {\bar A}$, $\hat{\bar A} = S/Q$. 
 
 \[
 \begin{array}{ccccccc}
    &   & A           &            & S                  & &\hole \\
    &   & \downarrow  &            & \downarrow         & &\hole \\
A/I & = & \overline A & \subseteq  & \hat{\overline A}  &=& S/Q   \\
 \end{array}
 \]

Let $\varphi'$ be a lifting of $\overline \varphi$ in $S$. As by assumption $J(\bar x\cdot\overline \varphi) = 0$ and $\text{det}(\overline \varphi) \neq 0$, we have $J(x'\cdot\varphi') \subseteq Q$ and $\text{det}(\varphi') \notin Q$.
 
Thus $\text{det}(\varphi') \notin J(x'\cdot\varphi')$, the sequence $x'\cdot\varphi'$ is a maximal regular sequence in $S$, see proposition \ref{p1.4}, and $\text{det}(\varphi')$ belongs to every ideal of $S$ containing strictly $J(x'\cdot\varphi')$, see corollary \ref{c1.2}(iv). But $\text{det}(\varphi') \notin Q$, so the inclusion $J(x'\cdot\varphi') \subseteq Q$ is in fact an equality, this takes care of (a). Moreover we just saw that $\hat {\bar A} = S/Q$ is a complete intersection of dimension 0, thus $\bar A$  also has dimension 0 and we have $ \bar A =\hat {\bar A}=  S/J(x'\cdot \varphi ')$. Now (b) follows from corollary \ref{c1.2}(i)(iii).

 $(iii) \Rightarrow (c).$ \quad  We already proved the equivalence of $(i),(ii), (iii)$. So let $\varphi$ be a matrix as in (iii). As $\text{det}(\varphi)\notin J(x\cdot\varphi)$, this equivalence, applied to the ideal $J(x\cdot\varphi)$ of A, tells us that this ideal is also a C.I.0-ideal. As we already proved $(iii)\Rightarrow (b)$, we  know that the socle of $A/J(x\cdot\varphi)$ is generated modulo $J(x\cdot\varphi)$ by $\text{det}(\varphi)$, which gives us the last equality in (c). But this socle is contained in every non null ideal of the zero-dimensional complete intersection $A/J(x\cdot\varphi)$. As $\text{det}(\varphi) \notin I$, the inclusion $J(x\cdot\varphi)\subseteq I$ is in fact an equality. 
\end{proof}

\begin{remark}\label{ci} Note that we have recovered, in dimension 0, a well-known fact: when a complete intersection is the quotient of a regular local ring $S$, it is the quotient of $S$ by an ideal generated by a regular sequence.
\end{remark}

We now retain the condition on $\overline\varphi$ met in the above proof.

\begin{definition}\label{d1.3}
Let $A$ be a noetherian local ring with a sequence $x = (x_1,\cdots,x_n)$ generating its maximal ideal (not necessarily minimally). An {\bf $x$-Wiebe matrix} for the ring $A$, if it exists, is an $x$-nice matrix $\psi$ belonging to the null ideal, i.e. a matrix $\psi$ such that $J(x\cdot\psi) = 0$ and $\text{det}(\psi) \neq 0$.
\end{definition}

We summarize the information given by the preceding theorem.

\begin{corollary}\label{crw}
A noetherian local ring $A$ is a complete intersection of dimension 0 if and only if it has an $x$-Wiebe matrix for some (every) sequence $x$ generating its maximal ideal. When this is the case, the determinant of an $x$-Wiebe matrix generates the socle of the ring.

More generally, for any noetherian local ring $A$ with a sequence $x$ generating its maximal ideal, an ideal $I$ of $A$ is a C.I.0-ideal of $A$ if and only if there is an $x$-nice matrix $\varphi$ belonging to it, i.e. a matrix $\varphi$ such that $det(\varphi) \notin J(x\cdot\varphi) = I$. When this is the case, let us denote by $\overline{(\cdot )}$ the images modulo $I$. The matrix $\overline{\varphi}$ is then an $\bar x$-Wiebe matrix for $\bar A=A/I$ and $socle(\bar A) = det(\overline {\varphi})\bar A$.
\end{corollary}

\begin{recall}\label{rg} In a zero-dimensional Gorenstein local ring $A$, the function which assigns to each ideal of $A$ its annihilator is an order-reversing involution of the set of ideals of $A$, see \ref{gor}. Under this involution,  the Gorenstein ideals (i.e. the proper ideals such that the corresponding quotient is Gorenstein ) correspond to the nonnull principal one's, see \cite{scac}. This in essence goes  back to Macaulay (\cite{MAC}, Ch. IV), see also (\cite{ms}, II.2 and VI.2).
\end{recall}

 As a C.I.0-ideal of $A$ is also a Gorenstein ideal, this raises a question.
  For which $b\in A$ is $\text{Ann}(b) = 0:bA$ a C.I.0-ideal ? 
  
  The question makes sense whatever the local ring $A$ is, though only local rings of depth null do have an element $b$ the annihilator of which is a C.I.0-ideal. The answer again lies in Wiebe's criterion. More information on  this will be provided in corollary \ref{cwa}.

\begin{corollary}\label{cann}
Let $A$ be a noetherian local ring with a sequence $x =(x_1,\cdots,x_n)$ generating its maximal ideal M, and let $0\neq b\in A$.The following conditions are equivalent:
\begin{enumerate}
\item[(i)] $Ann(b)$ is a C.I.0-ideal,
\item[(ii)]  $\exists  \varphi \in A^{n\times n}$ such that \quad $b\cdot x\cdot\varphi = (0,\cdots,0)$ \quad and \quad $b\cdot det(\varphi) \neq 0$,
\item[(iii)] $b\cdot \delta^A_0 (b\cdot M) \neq 0$.
\end{enumerate}
Moreover, any matrix $\varphi$ satisfying (ii) is an $x$-nice matrix belonging to $Ann(b)$: \quad $J(x\cdot\varphi) = 0:bA$ and $J(x\cdot \varphi):M=J(x\cdot \varphi)+det(\varphi)A$.
\end{corollary}
\begin{proof}
$(i)\Leftrightarrow (ii)$ and the last assertion. This follows from theorem \ref{twi}. Indeed, the condition ``$J(x\cdot\varphi) \subseteq \text{Ann}(b)$ and $\text{det}(\varphi) \notin \text{Ann}(b)$" is equivalent with (ii).

$(ii)\Leftrightarrow (iii).$ This follows from the description of the initial Fitting ideal given in \ref{rfit}.
\end{proof}

We now take some time for examples and remarks.

\begin{remark}\label{rwi}
To obtain the $x$-Wiebe matrices of the zero-dimensional complete intersection A with its sequence $x = (x_1,\cdots,x_n)$ generating its maximal ideal, we take an $x$-realization  $(S,Q,x')$ of $A$. With the preceding theorem, applied to the ideal $Q$ of $S$, we know that this ideal $Q$ is generated by a maximal regular sequence, say $a' = (a'_1,\cdots,a'_n)$. We can write $a' = x'\cdot\psi'$, for some matrix $\psi' \in S^{n\times n}$. The image of $\psi'$ in $A$ is an $x$-Wiebe matrix of $A$ and it follows from theorem \ref{twi} that every $x$-Wiebe matrix is obtained in that way.

\wl Thus an $x$-Wiebe matrix for the ring $A$ encodes the structure of $A$. This is particularly clear in the case where the ring $A$ is of equal characteristic and has thus a coefficient field $K$. In that case $A$ is a quotient of a ring of formal power series  $S=K[[x'_1, \cdots ,x_n']]$, $A=S/Q$ for an ideal $Q$ of $S$ generated by a maximal regular sequence, every element $a\in A$ can be written as a polynomial expression in the $x_i$'s with coefficients in $K$ and any such writing gives us an element $a'\in S'$ the image of which in $S/Q$ is $a$. Given an $x$-Wiebe matrix $\psi$ of $A$, it is then easy to lift it into an $x'$-nice matrix $\psi '\in S^{n\times n}$ and once this is done we have a minimal  set of generators $x'\cdot \psi '$ of the ideal $Q= J(x'\cdot \psi ')$. The most simple is the $x$-Wiebe matrix $\psi$, the most simple is our minimal set of generators of the ideal $Q$.

\wl We  note that, if $\psi$ is an $x$-Wiebe matrix for the ring $A$ 	and if $\theta$ is an invertible matrix of the same size, then $\psi\theta$ is also an $x$-Wiebe matrix for $A$. But the set of the  $x$-Wiebe matrices of $A$ may be larger than the set $\{ \psi\theta \mid \theta \quad \text{invertible}\}$.

(The equivalence class of an $x$-nice matrix will be described in \ref{c42} .) 

\wl Note also that the invertible matrices of the right size are exactly the $x$-nice matrices belonging to the maximal ideal of $A$.
\end{remark}

\begin{example}\label{e1.1}
Let $A = K[X,Y]/(XY, X^3 +Y^3)$, where $K$ is a field, and let us denote by $x$, $y$ the images of $X$ and $Y$ in $A$. This ring $A$ is local with maximal ideal generated by the sequence $(x,y)$ and is a complete intersection of dimension zero. Here are two $(x,y)$-Wiebe matrices for $A$
\[
\psi_1
=
\left(
\begin{array}{cc}
y & x^2\\
0 & y^2\\
\end{array}
\right)
\qquad
\psi_2
=
\left(
\begin{array}{cc}
0 & x^2\\
x & y^2\\
\end{array}
\right) .
\]
We do not have a matrix $\theta$ such that $\psi_2 = \psi_1\theta$ (look at the (2,1)-entries: $x \notin y^2 A$).
\end{example}

\begin{remark}\label{rnice}
In a noetherian local ring with a sequence $x$ generating its maximal ideal, we can have a matrix $\beta$ of the right size which is not $x$-nice, though the ideal $J(x\cdot\beta)$
 is a C.I.0-ideal. But then we shall have another matrix, say $\varphi$, not unique, which is $x$-nice and for which $J(x\cdot\beta) = J(x\cdot\varphi)$. The easiest example of this phenomenon is  the null matrix of the right size in a complete intersection of  dimension 0.
\end{remark}

\begin{example}\label{ewi}
(Wiebe) Let $A = K[[X,Y,Z]]/(X^2,Y^2,Z^2)$, where K is a field. As usual, denote by $x$, $y$, $z$ the images of $X$, $Y$, $Z$ in $A$. The ring $A$ is local with maximal ideal $M$ generated by the sequence $(x,y,z)$. 

Here the matrix $\left( \begin{array}{ccc} x & 0 & 0\\ 0 & y & 0\\ 0 & 0 & z\\ \end{array}\right)$ is an $(x,y,z)$-Wiebe matrix of $A$,

\wl and we have: $0:M = xyzA = M^3 \neq 0$, $M^4 = 0$.

\wl Let $c\in K$ and $\varphi = \left(\begin{array}{ccc} 1 & 0 & 0\\ -c & y & 0\\ 0 & 0 & z\\ \end{array} \right)$.

\wl We have: $(x+cy)\cdot (x,y,z)\cdot\varphi = (cxy,xy,xz+cyz)\cdot\varphi = (0,0,0)$ and $(x+cy)\cdot\det(\varphi)= xyz \neq 0$. With corollary \ref{cann}, we find that $\text{Ann}(x+cy) = J((x,y,z)\cdot\varphi)= (x-cy)A$ is a C.I.0-ideal.(Direct computations work as well.)

Look now at the element $x+y+z$ of $A$. When the characteristic of $K$ in not 2, a direct computation shows that $\text{Ann}(x+y+z) = J( xy-xz, xy-yz)$, thus the Gorenstein local ring $A/\text{Ann}(x+y+z)$ has embedding dimension 3, while the third power of its maximal ideal is null, such a ring cannot be a complete intersection, see corollary \ref{c1.3}.
\end{example}

We end this section with  small but funny applications of the ideas developed here.

\begin{proposition}\label{e1.2}
Let $A$ be a noetherian local ring of embedding dimension 2 and let $(x,y)$ be a sequence generating its maximal ideal. If $xy = 0$ and $xA\cap yA \neq 0$, then $A$ is a complete intersection of dimension 0.
\end{proposition}
\begin{proof}Indeed, $\exists c,d \in A$ with $xc=yd \neq 0$. Now the matrix
$
\left( 
\begin{array}{cc}
c & 0\\ 
{-d} & x \\ 
\end{array} 
\right)
$ 
is an $(x,y)$-Wiebe matrix of $A$.
\end{proof} 

Note that the conditions  $xy = 0$ and $xA\cap yA \neq 0$ are satisfied for the ring of the example \ref{e1.1}.
\vspace{0.5 cm}

In example \ref{ewi} we have seen  a complete intersection of dimension zero with a Gorenstein ideal which is not a C.I.0-ideal. Such a phenomenon cannot occur in embedding dimension two. Indeed, a grade two ideal $I$ of a regular local ring $S$ such that the quotient ring $S/I$ is Gorenstein is always a complete intersection ideal.

In particular, we have the following proposition, which again goes back to Macaulay. It will be useful to provide and understand some examples (see \ref{eouf}). We present a rather unusual proof of it based on the preceding considerations.

\begin{proposition}\label{pgor}
A Gorenstein local ring $A$ of dimension 0 and embedding dimension at most 2 is a complete intersection.
\end{proposition}
\begin{proof}
The case of embedding dimension less than 2 is obvious. So, let $(x,y)$ be a sequence generating minimally the maximal ideal $M$ of $A$. The quotient ring $A/xA$, having embedding dimension 1, is the quotient of a discrete valuation ring. Thus the zero-dimensional ring $A/xA$ is a complete intersection and Gorenstein, and the ideal $0:xA$ is principal, say $0:xA = zA$ for some $z\in M$.
Let $c\in A$ be an element which generates, modulo $xA$, the socle of $A/xA$.

We have: $xz = 0$, $xA = 0:zA$, $c\notin xA$, $c\cdot M \subseteq xA$. Thus $cz \neq 0$ and $yc = xd$ for some $d\in A$.

With this we see that the matrix $\left(\begin{array}{cc} z & -d\\0 & c\\ \end{array} \right)$ is an $(x,y)$-Wiebe matrix for $A$.

We also note that the entries of this matrix are not invertible.
\end{proof} 

With the above proof we recover a probably rather well-known fact, which can also be viewed as a fact concerning the intersection of two algebraic plane curves, maybe slightly generalized. 

\begin{corollary}\label{embg} Let $S$ be a regular local ring of dimension two, of maximal ideal $M'$.

For any parameter ideal $Q\subset M'^2$ of $S$, for all $x'\in M'\setminus M'^2$, there exist $z'\in M', g'\in M'^2$ such that the sequence $(x'z',g')$ generates $Q$.
\end{corollary}
\begin{proof} Let $(x', y')$ be a sequence in $S$ generating $M'$ minimally. The ring $A=S/Q$ is a complete intersection of dimension zero and embedding dimension two, and its maximal ideal $M$ is minimally generated by the sequence $(x, y)$, where $x$, $y$ denote the images of $x'$, $y'$ in $A$. We lift the  $(x,y)$-Wiebe matrix for $A$ $\psi=\left(\begin{array}{cc} z & -d\\0 & c\\ \end{array} \right)$ given by the above proof into a matrix $\psi ' =\left(\begin{array}{cc} z' & -d'\\0 & c'\\ \end{array} \right) \in S^{n\times n}$. With (\ref{twi} (a)) we have that the ideal $Q$ of $S$ is generated by the sequence $(x', y')\cdot \psi '= (x'z', c'y'-d'x')$.
\end{proof}

In higher dimension, we shall have an analogous result as soon as the quotient ring $S/Q$ has a minimal generator of its maximal ideal, the annihilator of which is C.I.0, see \ref{pgo}.

\section{Chains of Gorenstein ideals}  

In this intermediary section we collect some facts about chains of Gorenstein ideals. In an Gorenstein  ring of dimension zero they correspond to a factorization of an element $v$ generating the socle of the ring. As C.I.0-ideals are  Gorenstein, this will also yield a few basic facts about chains of C.I.0-ideals, the  object of the next section. And since in this paper we are mainly concerned with C.I.0-ideals, most of our examples will be C.I.0-examples.

\begin{observations}\label{og} Let $A$ be a zero-dimensional Gorenstein local ring .
 We recall  that the Gorenstein ideals of $A$ are exactly the annihilators of the nonnull principal one's, see  \ref{rg}, and that the maximal ideal $M$ is the annihilator of the socle  $0:M$ of $A$ which is principal. 

Thus an inclusion  of Gorenstein ideals in $A$ 
 
 $$I_0 \subseteq I_1$$
corresponds to an inclusion  of nonnull principal ideals

$$b_1A = 0:I_1 \quad \subseteq \quad 0:I_0 = b_0A$$

 In this situation, we write \quad \quad \quad \quad $b_1 = c \cdot b_0$ \quad \quad \quad for some $c$ in $A$. \\
 We then have  $$I_0:I_1 = I_0+cA.$$

 Indeed, $c \cdot I_1 \subset I_0$ because $b_0 \cdot c \cdot I_1 = b_1 \cdot I_1 = 0$. On the other hand, if $t \cdot I_1 \subset I_0$ for some $t \in A$, then $b_0 \cdot t \cdot I_1 = 0, \quad b_0 \cdot t \in 0:I_1 =  b_1A$ and $b_0\cdot t = b_1 \cdot a$ for some $a \in A$. But then we have $b_0 \cdot t =  b_0 \cdot c \cdot a$ and $t - ca \in I_0$. 
  
   And we conclude that the inclusion $I_0 \subseteq I_1$ is strict exactly when $c$ is not invertible ($I_0 = I_1 \Leftrightarrow I_0:I_1= A \Leftrightarrow I_0 + cA= A \Leftrightarrow c \quad \text {is invertible}$).
  
\wl Thus a factorization of any nonnull element $b\in A$ corresponds to a chain of Gorenstein ideals with $0:bA$ as its largest term and conversely.

\wl On the other hand, remember that a local artinian ring $A$ is Gorenstein of exponent $t+1$ if and only if $0:M=\text{socle}(A) = vA=M^t$ for some nonnull $v\in A$,  that this $v$ then is a product $v=a_t\cdots a_1$ with each $a_i \in M\setminus M^2$, and that this socle $vA$ is contained in every nonnull ideal of $A$. 

So let $v$ be an element of $A$ generating its socle.
For any Gorenstein ideal $I_0=0:b_0A$ of our Gorenstein local ring $A$, $b_0\neq 0$, we thus can write $v=hb_0$ for some $h\in A$. 
 By the above such an $h$ generates modulo $I_0$ the socle of $A/I_0$  and every factorization 
$$h=c_j\cdots c_1$$
corresponds now to a chain of Gorenstein ideals with $I_0=0:b_0A$ as its smaller term
$$I_0 \subset I_1\subset\cdots \subset I_j=M$$
where $I_i = 0:(c_i\cdots c_1b_0)A$ for $1\leq i \leq j$.

\wl In particular, a strict saturated chain  $$0 \subset I_1 \subset \cdots \subset I_s = M$$ of Gorenstein ideals corresponds to a strict saturated chain of nonnull principal ideals and to a factorization $$v = a_s \cdots a_1$$ of any element $v\in A$ generating the socle of $A$, in such a way that
    $$I_i = 0: (a_i \cdots a_1)A$$
     We note that the elements $a_i$ occurring in the factorization are not invertible because the chain we started with was strict, and we conclude that $M^s \neq 0$. This gives us a bound on the length of chains of Gorenstein ideals in $A$. We also note that these $a_i$'s cannot be written as a product of two non invertible elements because the chain was saturated.
    
    Now, if our Gorenstein local ring $A$ has exponent $t+1$, i.e. if  $0:M=vA=M^t$, since this nonnull element $v$ generating the socle of $A$ may be written as a product of $t$ elements of $M$,
  we have a chain of Gorenstein ideals of the maximal length $t$. However, we may have a strict saturated chain of C.I.0-ideals of length less than $t$, as we shall see in the examples below.      
\end{observations}
 
 \wl First we summarize part of our observations.
 
 \begin{proposition}\label{pgs}Let $A$ be a zero-dimensional Gorenstein local ring with maximal ideal $M$.  
 
 (i)  A nonnull Gorenstein ideal $(0:bA)$ contains strictly another nonnull Gorenstein ideal exactly when the element $b$ can be written as a product of two non invertible elements.
 
 (ii) A chain of Gorenstein ideals $0\subset I_1\subset \cdots \subset I_r =M$ corresponds to a factorization $v=a_r\cdots a_1$ of an element $v$ generating the socle $(0:M)$ of $A$ in such a way that $I_i=0:(a_i\cdots a_1)$.
 
 In this correspondence we have \quad $I_{i-1}:I_i=I_{i-1}+a_iA$\\ and \quad $I_i:M=I_i+(a_r\cdots a_{i+1})A$.
 
 The chain is strict if and only if the $a_i$'s are non invertible.
 
 (iii) Assume that the exponent of $A$ is $t+1$, which means $M^t\neq 0$ and $M^{t+1}=0$, then any strict chain of Gorenstein ideals has length at most $t$ and there is a chain of the right length $t$.
  
  In particular, if \quad $0\neq b \in M^i$, then $M^{t+1-i}\subseteq (0:bA)$, so that $\text{exponent}(A/(0:bA))\leq (t+1-i)$  and a chain of Gorenstein ideals with $(0:bA)$ as its smaller term has length at most $t-i$.
 \end{proposition}
  
In the above, we note that the inequality "$\text{exponent}(A/(0:bA))\leq (t+1-i)$" may be strict, see \ref{eiar}.

\begin{example}\label{etoc}  Let $A = K[X,Y,Z,T]/(X^2,Y^2,Z^2,T^2)$, where $K$ is a field . As usual, denote by $x$, $y$, $z$, $t$ the images of $X$, $Y$, $Z$, $T$ in $A$. This ring $A$ is a Gorenstein local ring (and even a complete intersection of dimension zero) with maximal ideal $M$ generated by the sequence $ (x,y,z,t)$. Here $M^4 \neq 0, M^5=0$ and $0:M=(xyzt)A$, our ring $A$ has exponent 5 and we have a strict chain of Gorenstein ideals of length 4.

But we shall provide a nonnull element of $M^2$ which cannot be written as the product of two non invertible elements of $A$. With the above proposition we know that the annihilator of such an element is a Gorenstein ideal, minimal among the nonnull Gorenstein ideals of $A$, and we also know that every strict saturated chain of Gorenstein ideals passing through this annihilator  has length at most 3.

 We claim that $(xy+xz+zt)$ is such an element. Indeed, as our ring $A$ is graded, if $(xy+xz+zt)$ is the product of two non invertible elements, it is also the product of two elements of degree 1 and we have

$xy+xz+zt=(ax+by+cz+dt)\cdot (a'x+b'y+c'z+d't),$\\
for some $a,b,c,d,a',b',c',d' \in K$.

We must have \quad $ab'+ba'=1$, so one at least of the two elements $ab',ba'$ is invertible, say $ab'\neq 0$.

 We must have \quad $cd'+dc'=1$, so again one of the two elements $cd',dc'$ is nonnull, say $cd'\neq 0$ (the case where $dc'$ is nonnull is similar to this one).
 
 We also have \quad $bc'+cb'=0$. Since $c\neq 0\neq b'$, we also have $bc'\neq 0$ and we may write \quad $b'=ub,\quad c'=-uc  \quad \text{for some}\quad 0\neq u \in K$.
 
 We now have \quad $0=bd'+db'=bd'+dub$ \quad and we already know that $b\neq 0\neq d'$. Thus we also have \quad $d'=-ud$ \quad and \quad $d\neq 0$.
 
  We also have \quad $ad'+da'=0$. With the preceding we then have \quad $-uad+a'd=0$ \quad and we already know $d\neq 0$. We obtain \quad $a'=ua$.

 Finally, since we also have  \quad $ac'+ca'=1$, we obtain \quad $-uac+ca'=1$ \quad and \quad $a'-ua=c^{-1}$,\quad in contradiction with \quad $a'=ua$

Note that the Gorenstein ideal $I=0:(xy+xz+zt)A$ is not a C.I.0-ideal. Indeed, computations show that $I\subset M^2$, thus the quotient ring $A/I$
 has embedding dimension 4, while its exponent is at most 4 since $0\neq I \supset \text{Socle}(A)=M^4$. Such a ring cannot be a complete intersection of dimension zero, see \ref{c1.3}. (In fact, with \ref{rgra}, we exactly have $\text{exponent}(A/I)=3$.) 
\end{example}

\begin{example}\label{eouf} Let $A=\mathbb Q[X,Y]/(X^3, Y^3]]$, where $\mathbb Q$ is the field of rational numbers. This ring is  Gorenstein local, (it is even a complete intersection of dimension zero), with our usual notations its maximal ideal $M$ is generated by the sequence $(x,y)$ and  we have $M^4 \neq 0,\quad M^5=0$ and $(0:M)=(x^2y^2)A$, this ring  $A$ has exponent 5. 

With arguments similar to those used in the preceding example, we see that the element $x^2+y^2$ cannot be written as the product of two elements in $M$. Here again, as in the preceding example, any strict saturated chain of Gorenstein ideals passing through $(0:(x^2+y^2)A)$ has length at most 3, though we have other strict chains of length 4.

But this example is also a C.I.0-example. Indeed,  the embedding dimension of $A$ is 2, so that every Gorenstein ideal of $A$ is also a C.I.0-ideal, see proposition \ref{pgor}. We have here a complete intersection of dimension zero and exponent 5 with a strict saturated chain of C.I.0-ideals of length less than 4.

Let us now provide an $(x,y)$-nice matrix belonging to $(0:(x^2+y^2)A)$.
Let $\varphi = \left( \begin {array}{cc}  y & -x\\0 & y \end {array} \right)$. We have 
$(x^2+y^2)\cdot (x,y)\cdot \varphi = (x^2+y^2)\cdot (xy, -x^2+y^2) =(0,0) $ and $(x^2+y^2) \text{det}(\varphi)=x^2y^2 \neq 0$. 
We conclude with \ref{cann} that $\varphi$ is an $x$-nice matrix belonging to $(0:(x^2+y^2)A) = (xy, y^2-x^2)$.

Finally, let us here exhibit another phenomenon. Under the correspondence between chains of Gorenstein ideals and factorizations of an element $v$ generating the socle observed in \ref{og}, it is quite clear that a strict saturated chain of Gorenstein ideals corresponds to a factorization of  $v$ which cannot be refined, in which all the factors are not products of two non invertible elements. However the converse is not true. In our ring $A$, the socle is generated by $x^2y^2$. The length three factorization $x^2y^2= (x^2+y^2)\cdot y\cdot y$ cannot be refined, though the corresponding chain $0 \subset (0:yA)\subset (0:y^2A) \subset (0:vA)=M$ is not saturated and may be refined to $0 \subset (0:yA)\subset (0:y^2A) \subset  (0:y^2xA)\subset (0:vA)=M$, this extended chain corresponding to the length four factorization $v=x\cdot x\cdot y\cdot y$.

This example will be revisited in \ref{estrc}.
\end{example}

\begin{remark}\label{rgra} It is worthwhile to note that the zero-dimensional Gorenstein local rings rings in the preceding examples are {\bf positively graded}, with the maximal ideal generated  by the homogeneous elements of degree 1, while their homogeneous elements of degree zero form a field . 

 A zero-dimensional Gorenstein local ring, when positively graded, enjoys more symmetry and has nicer properties than a non-graded one. We now recall some of them.
 
 So let $A$ be a positively graded zero-dimensional Gorenstein local ring,  and let $t+1$ be the exponent of $A$, which means that its maximal ideal satisfies $M^t\neq 0$, $M^{t+1}=0$.
 
We first note that the socle $(0:M)=M^t$ of $A$, which is principal isomorphic to $K=A/M$,  is a homogeneous ideal. We also note that, if $0\neq v\in M^t$, then $v$ is homogeneous of degree $t$ and generates $M^t$.
 
  Moreover, if we have a factorization of such a $v$ generating $M^t$,  $v=hb$, where  $b$ is homogeneous of degree $i, \quad 0 <  i\leq t$, then we also have $v=h_{t-i}b$, where $h_j$ denotes the degree $j$ component of $h$. Thus $h_{t-i}$ generates modulo $(0:bA)$ the socle of the graded Gorenstein local ring $A/(0:bA)$. This implies that $M^{t-i}\nsubseteq (0:bA)$ and that $ (A/(0:bA))$ has exponent exactly $(t+1-i)$ (since obviously  $M^{t+1-i}\subset (0:bA)$).

 With the above we also can see  that the multiplication in $A$ induces an exact pairing 
$$M^i/M^{i+1}\times M^{t-i}/M^{t-i+1} \rightarrow M^t\simeq K=A/M$$
so that the $K$-vector spaces $M^i/M^{i+1}$ and $ M^{t-i}/M^{t-i+1}$ have the same finite dimension. This again was already observed by Macaulay  in \cite{MAC}.

A non graded example without the above property will be given in \ref{eiar}.
\end{remark}

  We now provide some easy facts concerning chains of Gorenstein ideals of the maximal length and those minimal nonnull Gorenstein ideals which can be the starting point of a such a chain. This will also be useful in the next section where we shall
 look at the more particular chains of C.I.0-ideals.

\begin{proposition}\label{equi} Let $A$ be a Gorenstein local ring of dimension zero and exponent $t+1$ and let $v$ be an element generating the socle of $A$:\quad $0:M= M^t =vA$.

(i) \quad Let 
$0\subset I_1 \subset \cdots \subset I_t=M$ be a strict chain of Gorenstein ideals of the maximal length $t$
and let $v=d_t\cdots d_1$ be a factorization of  $v$  associated to this chain as in \ref {pgs}.

Then each $d_i$ is a minimal generator of the maximal ideal $M$ of $A$ and $\text{exponent}(A/I_i)=\text{exponent}(A)-i$.

 In particular the Gorenstein ideal  $I_1$ is the annihilator of a minimal generator of the maximal ideal   and  $\text{exponent}(A/I_1)=\text{exponent}(A)-1$.

(ii)(a) \quad Let $I$ be  another Gorenstein ideal of $A$ such that $\text{exponent}(A/I)=\text{exponent}(A)-1$.

 Then the principal ideal $(0:I)$ is generated by an element $y\in M\setminus M^2$.
 
 (b) \quad Conversely, if $A$ is positively graded with its maximal ideal generated by the homogeneous elements of degree 1, and if $y$ is a homogeneous minimal generator of the maximal ideal, then  the Gorenstein ideal $(0:yA)$ has the property that  $\text{exponent}(A/(0:yA))=\text{exponent}(A)-1$.
\end{proposition} 

\begin{proof} (i) \quad  With  \ref{og} or \ref{pgs} we know that the $d_i$'s are non invertible. Since  their product  generates $M^t$, we must have $d_i\notin M^2$\quad $\forall i, 1\leq i \leq t$,. 

 With \ref{og} or \ref{pgs} we have $(0:d_1A) =I_1$.

Since the nonnull ideal $I_1$ contains the socle $0:M=M^t$ of $A$, we have $\text{exponent} (A/I_1)\leq t$. On the other hand, $A/I_1$ contains a strict chain of Gorenstein ideals of length $t-1$, so with \ref{pgs} we also have $\text{exponent}(A/I_1)\geq t$. The general assertion on $\text{exponent}(A/I_i)$ is obtained by iteration.  

(ii)(a)  Let $(0:I)= yA$ and write $v=hy$, where $v$ generates the socle $M^t$ of $A$. By hypothesis we have $M^{t-1} \nsubseteq I$ and we know with \ref{og} that $h$ generates modulo $I$ the socle  of $A/I$ which is contained in any nonnull ideal of $A/I$. From this we conclude $h \in M^{t-1}+I$.
This, together with $yI=0$ and $v=hy$, implies $y\notin M^2$, but $y\in M$ since $I\neq 0$. 

(b) This follows from \ref{rgra} since here $y$ is homogeneous of degree 1.
\end{proof}

As announced after \ref{pgs}, we note that the converse of (ii)(a) in the above proposition fails in general, see the following.

\begin{example}\label{eiar} Let $A=  K[X,Y]/(XY, Y^2- X^3)$, where $K$ is a field and let $x,y$ be as usual the images of $X$ and $Y$ in A.

This ring is again a complete intersection of dimension zero, thus Gorenstein. With our usual notations its maximal ideal $M$ is generated by the sequence $(x,y)$ and we have $xy=0$, $y^2=x^3$, $y^3=x^4=0$, $\text{socle}(A) = M^3 = x^3A= y^2A$. This ring  has exponent 4. 

In this ring we have $y\cdot (x, y)\cdot \left(\begin{array}{cc}1&0\\0&y \end{array}\right)= (0,0)$ and $y\cdot \text{det}\left(\begin{array}{cc}1&0\\0&y \end{array}\right)\neq 0$. With \ref{cann}, this tells us that the Gorenstein ideal $0:yA$ is also a C.I.0-ideal, that $0:yA= J((x, y)\cdot \left(\begin{array}{cc}1&0\\0&y \end{array}\right))=xA$. But
 the quotient ring $A/xA \simeq K[Y]/(Y^2)$ has exponent $2<4-1$, though $y\in M\setminus M^2$.

\end{example}

 \section{Matrix factorizations and chains of C.I.0-ideals} 
 
  Now  we turn to chains of C.I.0-ideals in any noetherian local ring $A$ with a sequence $x = (x_1, \cdots , x_n)$ generating its maximal ideal. We know that  C.I.0-ideals correspond to $x$-nice matrices, though this correspondence is not bijective, there are a lot of $x$-nice matrices belonging to the same C.I.0-ideal.  Anyway, we shall see that the factorizations of the $x$-nice matrices belonging to a C.I.0-ideal $I$  correspond to the ascending chains of C.I.0-ideals starting at $I$,  just like the factorizations of an element $v$ generating the socle of a zero-dimensional Gorenstein local ring $A$ correspond to the chains of Gorenstein ideals in $A$.

\wl  First we look at an inclusion $I_0\subseteq I_1$ of C.I.0-ideals in $A$ and we explore the relationships between the $x$-nice matrices belonging to $I_1$ and those belonging to $I_0$.

\begin{proposition} \label{pfac}
Let $A$ be a noetherian local ring with a sequence $x = (x_1, \dots ,x_n)$ generating its maximal ideal $M$. Let $\varphi_0$ be an $x$-nice matrix belonging to a C.I.0-ideal $I_0$  and assume we have a factorization $\varphi_0 = \varphi_1 \cdot \gamma$.

Then  we have:

(i) $\varphi_1$ is an $x$-nice matrix and $J(x\cdot \varphi_0)\subseteq J(x\cdot \varphi_1)$,

(ii) $J(x\cdot \varphi_1):M = J(x\cdot \varphi_1) + det(\varphi_1)A$,

(iii) $J(x\cdot \varphi_0) : det(\gamma) = J(x\cdot \varphi_1)$,

(iv) $J(x\cdot \varphi_0) : J(x\cdot \varphi_1) = J(x\cdot \varphi_0) + det(\gamma)A$.

Moreover, the inclusion (i) is strict if and only if the matrix $\gamma$ is not invertible, while $J(x\cdot \varphi _1)\neq M$ if and only if $\varphi _1$ is not invertible.
\end{proposition}
\begin{proof}
Assume that the matrix $\varphi_1$ is not $x$-nice. Then $\text{det}(\varphi_1) \in J(x\cdot \varphi_1)$ and there is a column matrix $\alpha \in A^{n\times 1}$ with $\text{det}(\varphi_1) = x\cdot \varphi_1 \cdot \alpha$. We then have 
\[\text{det}(\varphi_0) = \text{det}(\gamma)\cdot  \text{det}(\varphi_1) = \text{det}(\gamma)\cdot x \cdot \varphi_1 \cdot \alpha = (x\cdot \varphi_1 \cdot \gamma \cdot \gamma^c)\cdot \alpha = (x\cdot \varphi_0)\cdot (\gamma^c \cdot \alpha)
\]
where $\gamma^c$ denotes the adjoint matrix of $\gamma$.
This means that $\text{det}(\varphi_0) \in J(x\cdot \varphi_0)$, in contradiction with the hypothesis on $\varphi_0$. Thus $\varphi_1$ is x-nice and $I_1=J(x\cdot \varphi _1)$ is a C.I.0-ideal containing $I_0$.  

The equality (ii) is already known, see corollary \ref{crw}.

We also have: $\text{det}(\gamma)\cdot x \cdot \varphi_1 = x \cdot \varphi_1 \cdot \gamma \cdot \gamma^c = (x \cdot \varphi_0) \cdot \gamma^c$, so that $\text{det}(\gamma)\cdot J( x \cdot \varphi_1)\subseteq J(x\cdot \varphi_0)$. Moreover $\text{det}(\gamma)\cdot \text{det}(\varphi_1) = \text{det}(\varphi_0) \notin J(x\cdot \varphi_0)$. Corollary \ref{cann} applied to the ring $A/I_0$ and the image of $\text{det}(\gamma)$ in that ring gives the third equality, and the last one follows with recall \ref{rg}

Finally, when $\gamma$ is not invertible we see with equality (iv) that the inclusion (i) is strict and when $\varphi _1$ is not invertible we see with equality (ii) that $J(x\cdot \varphi _1)\neq M$.
\end{proof}
 As we are looking at matrix factorizations, the following terminology will be convenient.
 
  \begin{definition}\label{dind}Let $A$ be a commutative ring (with unit) and let $\gamma \in A^{n\times n}$ be a non invertible matrix. We say that this matrix $\gamma$ is {\bf decomposable} (in $A$) if it can be written as the product of two non invertible matrices of $A^{n\times n}$. If not, we say that $\gamma$ is {\bf indecomposable} (in $A$).
  \end{definition}
  
 \begin{remark}\label{rwstr} The above proposition tells us that, when an $x$-nice matrix $\varphi _0$ is decomposable, the corresponding ideal $I_0$ is strictly contained in another C.I.0-ideal $I_1\neq M$. We shall see in a moment that, if we have a strict chain of C.I.0-ideals $I_0\varsubsetneqq I_1 \varsubsetneqq M$, then there is a decomposable  $x$-nice matrix belonging to $I_0$. However some other $x$-nice matrices belonging to this same ideal $I_0$ may be indecomposable. The easiest example of this phenomenon is an indecomposable $x$-Wiebe matrix for a local ring of embedding dimension at least two. ( Observe that a noetherian local ring of embedding dimension at least two always has a non trivial quotient of embedding dimension one and thus a non trivial C.I.0-ideal $I$, $0\neq I \neq M$, since a noetherian local ring of embedding dimension one is necessarily a complete intersection). Such an example will be given at the end of this section, see \ref{ewiind}.
 \end{remark}
  
\wl Before going the other way, we need some preliminaries.

\begin{lemma}\label{l} Let $A$ be a noetherian local ring with a sequence $x = (x_1, \cdots ,x_n)$ generating its maximal ideal, let $I$ be an ideal of $A$ and let us denote by $\overline{(\cdot)}$ the images modulo I.

Let $\varphi \in A^{n\times n}$ be a matrix such that the matrix $\overline{\varphi}$ is an $\overline{x}$-nice matrix. 

Then the matrix $\varphi$ itself is an $x$-nice matrix and $J(x\cdot \varphi) \supseteq I$.
\end{lemma}

\begin{proof}
The hypothesis on $\overline{\varphi}$ means that $\text{det}(\varphi) \notin J(x\cdot \varphi) + I$. We apply theorem \ref{twi}((iii), (c)) to the ideal $J(x\cdot \varphi) + I$ of $A$ and obtain that the matrix $\varphi$ is an $x$-nice matrix belonging to the C.I.0-ideal $J(x\cdot \varphi) + I$,  that $J(x\cdot \varphi) = J(x\cdot \varphi) + I$, so that $I \subseteq J(x\cdot \varphi)$.
\end{proof}

\begin{proposition}\label{pch} Let $A$ be a noetherian local ring with a sequence $x = (x_1,\cdots ,x_n)$ generating its maximal ideal $M$. Let $I_0 \subseteq I_1$ be two C.I.0-ideals of $A$ and let $\varphi_1$ be an $x$-nice matrix belonging to $I_1$, $I_1 = J(x\cdot \varphi_1)$.

Then there is a matrix $\gamma \in A^{n\times n}$ such that the matrix $\varphi_0 =\varphi_1 \cdot \gamma$ is an $x$-nice matrix belonging to $I_0$, $I_0 = J(x\cdot \varphi_0)$.
\end{proposition}

\begin{proof}
We denote by $\overline {(.)}$ the images modulo $I_0$ and we take an $\bar x$-realization  $(S,Q,x')$ of $\bar A$.

\[ \xymatrix {
A \ar[dr] &        & S\ar@{.>}[dl]\\
 A/I_0  \ar@{=}[r]         & \bar A \ar@{=}[r] &  S/Q
} \]

 With remark \ref{ci}, we know that the ideal $Q$ of $S$ is generated by a maximal regular sequence, say $b' = (b'_1, \cdots , b'_n)$.

We now take a matrix  $\varphi'_1 \in S^{n\times n}$ such that its image in $S/Q = \bar A$ is  $\overline{\varphi} _1$. As the matrix $\overline{\varphi}_1$ is an $\bar x$-matrix belonging to $\bar{I_1}$, the matrix $\varphi'_1$ is an $x'$-nice matrix and $J(x'\cdot \varphi'_1) \supseteq Q$, see the above lemma \ref{l}. We thus have a matrix   $\gamma' \in S^{n\times n}$ such that $b' = x'\cdot \varphi'_1 \cdot \gamma'$. As $\text{det}(\varphi'_1\cdot \gamma') \notin J(b') = Q$, see proposition \ref{p1.1}, the matrix $\varphi'_1 \cdot \gamma'$ is an $ x'$-nice matrix belonging to $Q$.

 Let now $\gamma \in A^{n\times n}$ be a matrix such that $\bar \gamma$ is the image of $\gamma'$ in $\bar {A} ^{n \times n}$, the matrix $\gamma$ is the one we are looking for. Indeed ,
 $\overline{\varphi}_1 \cdot \bar{\gamma}$ is an $\bar x$-Wiebe matrix for $\bar A$, thus $J(x\cdot \varphi_1 \cdot \gamma) \subset I_0$. On the other hand, with lemma \ref{l} again, we have that $\varphi_1 \cdot \gamma$ is $x$-nice and $I_0 \subset J(x\cdot \varphi_1 \cdot \gamma)$. Thus the matrix $\varphi_0 = \varphi_1 \cdot \gamma$ is an $x$-nice matrix belonging to $I_0$.
\end{proof}

The particular case when the ring $A$ itself is a complete intersection of dimension zero deserves attention. In that case, when one puts $I_0=(0)$, propositions \ref{pfac} and \ref{pch} provide new informations about the C.I.0-ideals and their annihilators.

\begin{corollary}\label{cwa}
Let $A$ be a zero-dimensional complete intersection with a sequence $x = (x_1, \cdots , x_n)$ generating its maximal ideal $M$, not necessarily minimally. 

 (i). An ideal $I$ of $A$ is a C.I.0-ideal if and only if there are matrices $\varphi, \gamma \in A^{n\times n}$ such that $I = J(x\cdot \varphi)$ and $\varphi \cdot \gamma$ is an $x$-Wiebe matrix for $A$.

In that case, $\varphi$ is an $x$-nice matrix belonging to $I$,\quad $I:M = I + det( \varphi)A$  and \quad $0:I = det(\gamma)A$.

 (ii). The annihilator of an element $b\in A$ is a C.I.0-ideal if and only if there are matrices $\varphi, \gamma \in A^{n\times n}$ such that $b = det(\gamma)$ and  $\varphi \cdot \gamma$ is an $x$-Wiebe matrix for $A$.

In that case, the matrix $\varphi$ is $x$-nice and $J(x\cdot \varphi) = 0:bA$.
\end{corollary}

We now look at chains of C.I.0-ideals. With  \ref{pfac} and \ref{pch}, we obtain a C.I.0 version of \ref{pgs}.

\begin{corollary}\label{cch}
Let $A$ be a zero-dimensional complete intersection with a sequence $x = (x_1, \cdots , x_n)$ generating its maximal ideal $M$, not necessarily minimally.

To every factorization in $A^{n \times n}$ of some x-Wiebe matrix of $A$
$$\psi = \eta_{t} \cdots \eta_1$$
corresponds a chain of C.I.0-ideals
 $$0= I_0\subset I_1 \subset \cdots \subset I_{t-1} \subset I_t = M$$ where \quad $I_i = J(x\cdot \eta_{t} \cdots \eta_{i+1})= 0:\text{det}(\eta_i\cdots \eta_1)$ for $1\leq i < t$.
 
 This correspondence is onto, but far from one to one, nevertheless the above chain is strict if and only if the matrices $\eta_i$ are not invertible for $1\leq i \leq t$.
 
 Under this correspondence, we also have $$I_{i-1}:I_i=I_{i-1}+\text{det}(\eta _i)A  \quad \text{and} \quad I_i:M=I_i+\text{det}(\eta _t)\cdots \text{det}(\eta _{i+1})A.$$
 \end{corollary}

\begin{example}\label{beiar} Let $A=  K[X,Y]/(XY, Y^2- X^3)$ be the ring of \ref{eiar}.

Here is a factorization of an $(x,y)$-Wiebe matrix corresponding as in \ref{cch} to the chain $0\subset xA\subset M$:
$$\left(\begin{array}{cc} y&-x^2\\0&y \end{array}\right ) =\left(\begin{array}{cc} 1&-x^2\\0&y \end{array}\right ) \cdot \left(\begin{array}{cc} y&0\\0&1 \end{array}\right )$$
\end{example}
 
 \begin{remark}\label{rcst}In the situation of the corollary \ref{cch}, given a factorization $\psi = \eta_{t} \cdots \eta_1$ of an $x$-Wiebe matrix $\psi$ of the ring, it is not easy to recognize at first sight if the corresponding chain of C.I.0-ideals is saturated or not. Of course, when this chain is saturated, all the matrices $\eta _i$ occurring in the factorization are indecomposable. However, as in the Gorenstein case (see \ref{eouf}), the converse is not true, it might happen that all the matrices $\eta _i$ are indecomposable while the corresponding chain is not saturated.
 
 Here again, as in \ref{rwstr}, the easiest example of this phenomenon is an indecomposable $x$-Wiebe in embedding dimension at least two. A less trivial example will be given at the end of this section (see \ref{estrc}).
 \end{remark}

Now we are faced with some problems concerning the length of the chains of C.I.0-ideals (which are also chains of Gorenstein ideals).

Recall that, in a complete intersection of dimension zero and exponent $t+1$, which means that the maximal ideal satisfies  $M^t \neq 0$ and $M^{t+1} = 0$, every strict chain of C.I.0-ideals has length at most $t$, see \ref{pgs} , and that we already have an example with a strict saturated chain of C.I.0-ideals of length less than $t$, see \ref{eouf}.
However, some chain problems remain, at least for us.

\begin{questions}\label{q3} Let $A$ be a complete intersection of dimension zero and exponent $t+1$, which means that its maximal ideal $M$ satisfies  $M^t \neq 0$ and $M^{t+1} = 0$.

 (i) \quad Do we always have a strict chain of C.I.0-ideals of the maximal length $t$~? If not, under which conditions do we have such a chain~?
 
(ii) \quad If moreover the embedding dimension of $A$ is $n$  (so that $n\leq t$, see \ref{c1.3}), do we have at least a strict chain of C.I.0-ideals of  length $n$~? If not, under which conditions~?

(iii) \quad   If moreover $\text{exponent}(A) = \text{embedding dimension}(A)+1 $, what can we say about these chains  ~?  
 \end{questions}

\begin{remarks}\label{rin} In the particular case when the complete intersection of dimension zero with a sequence $x$ generating its maximal ideal $M$ has a diagonal $x$-Wiebe matrix, all the entries of which are monomials in the $x_i$'s (such a ring could be called a monomial complete intersection of dimension zero) or just products of elements in $M\setminus M^2$, this matrix has a factorization of the right length and with \ref{cch} we do have a strict chain of C.I.0-ideals of the maximal length. Of course, the existence of such an $x$-Wiebe matrix has strong consequences for an $x$-realization $(S,Q,x')$ of our ring. With \ref{twi} it means that the ideal $Q$ of the regular local ring $S$ may be generated by elements which are completely factorizable, i.e. products of elements in $M'\setminus M'^2$ (where $M'$ denotes the maximal ideal of $S$).

The case of embedding dimension 1 is obvious because a local ring of embedding dimension 1 is a quotient of a discrete valuation ring. In that case every ideal is principal and C.I.0, there is only one strict saturated chain of C.I.0-ideals and it has the maximal length.

In embedding dimension 2, where Gorenstein ideals are C.I.0 (\ref{pgor}), we already know that there is a strict chain of C.I.0-ideals of the maximal length (\ref{pgs}). 

We already noticed that it is  easy to provide a non trivial C.I.0-ideal when the embedding dimension of our ring is at least 2.
 (Indeed, let  $x = (x_1, \cdots , x_n)$ be a sequence generating the maximal ideal $M$ minimally, with $ n\geq 2$ . The ideal $I$ generated by the sequence $(x_2, \cdots, x_n)$ is a non trivial ideal, $0\neq I\neq M$, it is also a C.I.0-ideal since the quotient $A/I$ has embedding dimension 1.)
 
Since C.I.0-ideals are Gorenstein, we already obtained  in \ref{equi} some information about the chains of C.I.0-ideals which have maximal length (when these exist) and those ideals which can be taken as the starting point of such a chain. In particular, \ref{equi} gives the following.
 \end{remarks}

\begin{proposition}\label{bequi} Let $A$ be a complete intersection of dimension zero and exponent $t+1$ and assume there is in $A$ a strict chain of C.I.0-ideals of the maximal length $t$.

Then there is in $A$ a minimal generator  of the maximal ideal $M$
 the annihilator of which is a C.I.0-ideal.

\end{proposition} 

\begin{proof} If \quad $0\subset I_1 \subset \cdots \subset I_t=M$ \quad is a chain of C.I.0-ideals and if 
 $\psi=\eta_t\cdots \eta_1$ is a factorization of an $x$-Wiebe matrix $\psi$ of $A$ associated to this chain as in \ref {cch}, then $I_1 =(0:\text{det}(\eta_1)A$. As  the factorization $\text{det}(\psi) = \text{det}(\eta _t) \cdots  \text{det}(\eta _1)$ is associated to the chain as in \ref{pgs}, we  then know with \ref{equi} that $\text{det}(\eta_1) \in M\setminus M^2$.
\end{proof}

Thus in \ref{q3} question (i) is related to another question concerning the minimal generators of the maximal ideal.
 
\begin{question}\label{q4} In a complete intersection $A$ of dimension zero, do we have a minimal generator $y$ of the maximal ideal, the annihilator of which is not only a Gorenstein ideal but also a C.I.0-ideal ? And what can we say about such a $y$, about the ring $A$, when  such a $y$ exists in $A$ ?
\end{question} 

These questions will be the object of the next section.

Now we just recall that we already have examples where the annihilator of some minimal generator of the maximal ideal is not a C.I.0-ideal, see Wiebe's example \ref{ewi}.

\wl Here are the other promised examples.

\begin{example}\label{ewiind} Let $A=\mathbb Q[X,Y]/(X^2,Y^2)$, where $\mathbb Q$ is the field of rational numbers. This ring is a complete intersection of dimension zero and exponent 3, with our usual notations its maximal ideal $M$ is generated by the sequence $(x, y)$ and $(0:M)=xyA$.

Here are two $(x,y)$-Wiebe matrices for $A$:

 $$\psi ^*= \begin{pmatrix}
x&0\\0&y
\end{pmatrix} \quad \text{and} \quad \psi= \begin{pmatrix}
x&-y\\y&x+y
\end{pmatrix} .$$

The first one is clearly decomposable, while the second one is not.

Let us see this. Assume $\psi$ is the product of two non-invertible square matrices, then both factors have determinants in $M\setminus M^2$ and , with lemma \ref{ltec}, we also have a factorization of one of the following forms
$$\psi = \beta \begin{pmatrix}
1&c\\0&d
\end{pmatrix} \quad \text{or} \quad \psi = \tilde\beta \begin{pmatrix}
e&0\\f&1
\end{pmatrix},\quad \text{where}\quad d, e\in M\setminus M^2.$$

In the first case, we have $\psi \begin{pmatrix}
1&-c\\0&1
\end{pmatrix} = \beta \begin{pmatrix}
1&0\\0&d
\end{pmatrix}$, which gives 
$$\left\{ \begin{array}{ccc}
-cx-y&=&ad\\
-cy+x+y&=&bd\\
\end{array}\right. ,\quad\text{where}\quad a=\beta_{12}, b=\beta_{22}.$$
As our ring is graded, we may take the degree one component of these equalities, we obtain
$$\left\{ \begin{array}{ccc}
-c_0x-y&=&a_0d_1\\
-c_0y+x+y&=&b_0d_1\\
\end{array}\right. ,\quad\text{where}\quad c_0, a_0, b_0 \in \mathbb Q.$$
From this we deduce $a_0\neq 0 \neq b_0$, putting $u=b_0^{-1}a_0$ we now obtain 
$$-c_0x-y=u(-c_0y+x+y),\quad \text{so that} \quad \left\{ \begin{array}{ccc}
-c_0&=&u\\
-1&=&u(-c_0+1)\\
\end{array}\right.$$
This gives us: \quad $u(u+1)=-1$. Since this last equation has no solutions in $\mathbb Q$, the first case is excluded.

\wl In the second case we have $\psi \begin{pmatrix}
1&0\\-f&1
\end{pmatrix} = \tilde \beta \begin{pmatrix}
e&0\\0&1
\end{pmatrix}$, so that $$\left\{ \begin{array}{ccc}
x+fy&=&a'e\\
y-f(x+y)&=&b'e\\
\end{array}\right. ,\quad\text{for some}\quad a', b' \in A.$$
As in the first case, we take the degree one components of these equalities, put $u={b'}_0^{-1}{a'}_0$ and   obtain
$$x+f_0y=u(y-f_0x-f_0y)\quad \text{ so that}  \quad \left\{\begin{array}{ccc}
1&=&-uf_0\\
f_0&=&u(1-f_0)
\end{array}
\right.$$ This again gives equations \quad $-u^{-1}=u(1+u^{-1})$, $-1=u(u+1)$ without solutions in $\mathbb Q$. The second case is also excluded.
\end{example}

\begin{example}\label{estrc} Let $A= \mathbb Q[X,Y]/(X^3, Y^3)$ be the ring of \ref{eouf}.

The matrix $\psi =\begin{pmatrix} x^2-y^2&-2xy\\xy&y^2+2x^2\end{pmatrix}$ is an $(x,y)$-Wiebe matrix for $A$ and has a factorization $$\psi=\eta_2 \cdot \eta_1,\quad \text{ where } \quad \eta_2= \begin{pmatrix}x&-y\\-y&2x\end{pmatrix}, \quad \eta_1=\begin{pmatrix}x&-y\\y&x\end{pmatrix}.$$
The elements $\text{det}(\eta_2)=2x^2-y^2$, $\text{det}(\eta_1)=x^2+y^2$ are both indecomposable (this can be seen with arguments similar to those used in \ref{etoc}, \ref{eouf}, we are working over the field of rational numbers). Thus the matrices $\eta_2$ and $\eta_1$ are also indecomposable. However, the corresponding chain 
$$0\subset (0:\text{det}(\eta_1)) \subset M$$
is not saturated. 

Indeed, let us put $I_1=(0:\text{det}(\eta_1))$. We already saw in \ref{eouf} that $I_1=(x^2-y^2, xy)$ is a C.I.0-ideal minimal among the nonnull C.I.0-ideals of $A$, but the ring $A/I_1 \simeq \mathbb Q[X,Y[/(X^2-Y^2, XY)$ has a non trivial C.I.0-ideal since it has embedding dimension 2.

Here is another factorization of another $(x,y)$-Wiebe matrix corresponding to a strict saturated chain of C.I.0-ideals starting at $I_1$
$$\tilde{\psi} = \begin{pmatrix}x&0\\-y&1\end{pmatrix} \cdot \begin{pmatrix}1&0\\0&x \end{pmatrix} \cdot \eta_1 = \begin{pmatrix}x^2&-xy\\0&x^2+y^2 \end{pmatrix}.$$

\end{example}

\section{About chains of C.I.0-ideals of the maximal length}

We have seen in \ref{bequi} that, in a complete intersection $A$ of dimension zero and exponent $t+1$, the existence of a chain of C.I.0-ideals of the maximal length t is related to the existence of some minimal generator $y$ of the maximal ideal $M$ of $A$, the annihilator of which is a C.I.0-ideal.

 So we first look at those minimal generators $y$ of the maximal ideal which have the above property, when they exist. It will turns out that the principal ideal $yA$ ($y\in M\setminus M^2)$ is itself C.I.0 as soon as its annihilator is C.I.0, that this condition on $y$ admits different formulations and has some consequences on the realizations of the ring $A$, see \ref{pgo}. This will be a consequence of our matrix factorizations together with a nice result of Kunz \cite{ku} .

\begin{theorem}\label{tku} (Kunz) Almost complete intersections are not Gorenstein.
\end{theorem}

In this theorem, a noetherian local ring is called an almost complete intersection if its completion with respect to the maximal adic topology is the quotient of a regular local ring by a grade $g$ ideal minimally generated by $g+1$ elements, i.e. an ideal minimally generated by a non regular sequence $(z_1, \cdots ,z_{g+1})$, where the subsequence $(z_1,\cdots ,z_g)$ is regular.  

\begin{corollary}\label{cku}In a complete intersection of dimension zero, a principal Gorenstein ideal is always a C.I.0-ideal and the annihilator of a principal C.I.0-ideal is again a principal C.I.0-ideal. 
\end{corollary}

When we factorize matrices, we sometimes encounter diagonal matrices. We need a notation for them.

\begin{notation}\label{not} An $s\times s$ diagonal matrix $\gamma$ with $\gamma _{ii}= a_i$ will be denoted by $\text{diag}(a_1, \cdots , a_n)$.
\end{notation}

To continue our program, we need an elementary technical lemma.  

\begin{lemma}\label{ltec} Let $A$ be a noetherian local ring with maximal ideal $M$ and let $\gamma \in A^{n\times n}$ be a matrix such that $d=\text{det}(\gamma)\in M\setminus M^2$.

(i)  There is an invertible matrix $\theta \in A^{n\times n}$such that the matrix $(\theta \gamma -I_n)$ has only one nonnull column (here $I_n$ denotes the $(n\times n)$-identity matrix), in other words such that $\theta \gamma$ has the following form $$ \left( \begin{array}{c|c|c}
I_r&\cdot &0\\
\hline
0&d&0\\
\hline
0&\cdot &I_s 

\end{array}
\right)
 , \quad r+s+1=n .$$
 
 (ii) There are invertible matrices $\theta_1, \theta_2 \in A^{n\times n}$ such that the matrix $\theta_1 \gamma \theta_2$ is diagonal of the form $\text{diag}(d, 1, \cdots ,1)$.
\end{lemma}

\begin{proof} We shall perform the usual operations on the rows (on the columns) of our matrix $\gamma$ which amount to multiply it on the left (on the right) by an invertible matrix.

(i)\quad First suppose that  all the entries $\gamma _{i1}$ in the first column of $\gamma$ belong to $M$, then at least one  of the corresponding $(n-1)\times (n-1)$ minors is invertible. After a row's permutation the minor belonging to the new $(1,1)$-entry is invertible. Then some operations on the last $n-1$ rows of this new matrix gives us a matrix of the form $\left( \begin{array}{c|c}
\cdot & \cdots\\
\hline
\vdots &1_{n-1}
\end{array} \right)$ One last row operation gives then the form described in (i) (with $r=0$).

This case having been handled, now assume that  one of the entries $\gamma _{i1}$ in the first column is invertible, a row operation bring it at the (1,1) place, some row's operations gives us a new matrix of the form $\left( \begin{array}{c|c}
1&\cdots \\
\hline
0&\tilde {\gamma}
\end{array} \right)$, where $\tilde{\gamma} \in A^{(n-1)\times (n-1)}$ and where $\text{det}(\gamma) =ud$ for some $u\notin M$. An induction on the size $n$ of the matrix allows us to bring the matrix $\tilde {\gamma}$ in the form described in (i), we then finish with some last row operations.

(ii) \quad Some column operations transform the matrix $\theta \gamma$ obtained in (i) into a diagonal matrix. We then bring its non-invertible entry at the first place with a row and a column operation.
\end{proof}

\begin{theorem}\label{pgo} Let $A$ be a complete intersection of dimension zero and let $y$ be any minimal generator of the maximal ideal $M$ of $A$. Let us put $I_1=(0:yA)$.

The following conditions for $y$ are equivalent:

(i) The ideal $I_1=(0:yA)$ is a C.I.0-ideal.

(ii) The ideal $I_1=(0:yA)$ is principal.

(iii)The principal ideal $yA$ is a C.I.0-ideal.

(iv) For any (for some) sequence $x=(x_1, \cdots ,x_n)$ generating the maximal ideal $M$ of $A$, there is an $x$-Wiebe matrix $\psi$ of $A$ such that its entries in the first column are all multiple of $y$, in other words there is an $x$-Wiebe matrix of the form $\psi=\varphi_1\cdot \text{diag}(y,1,\cdots ,1)$ (where $\varphi _1$ is an $x$-nice matrix belonging to $I_1$).

(v) For any (for some) sequence $x=(x_1, \cdots ,x_n)$ generating the maximal ideal $M$ of $A$, for any (for some) $x$-realization $(S,Q,x')$ of $A$, $A=S/Q$, and for all $y' \in S$ mapping onto $y$, the ideal $Q$ of $S$ is generated by a maximal regular sequence of the form $(y'z', a'_2, \cdots ,a'_n)$. 
\end{theorem}

\begin{proof} 

$(i) \Rightarrow (iv)$ \quad Let $x=(x_1, \cdots ,x_n)$ be any sequence generating the maximal ideal of $A$ and let $\varphi$ be an $x$-nice matrix belonging to the C.I.0-ideal $I_1=(0:yA)$. We have a matrix $\gamma \in A^{n\times n}$ such that $\tilde\psi = \varphi \gamma$ is an $x$-Wiebe matrix for $A$ and $I_1=(0:\text{det}(\gamma))$, see \ref{pch}, \ref{cwa}.  With \ref{rg} we then have  $yA=\text{det}(\gamma)A$ so that $y= \text{det}(\gamma)u$ for some invertible $u\in A$ (remember $y\in M \setminus M^2$). Thus $\text{det}(\gamma) \in M\setminus M^2$.

 With the lemma we have invertible matrices $\theta _1, \theta _2 \in A^{n\times n}$ such that $\theta _1 \gamma \theta _2 =\text{diag}(y, 1,\cdots , 1)$. Now the matrix $\psi=\tilde \psi \theta _2$ is still an $x$-Wiebe matrix for $A$ and has a factorization $\psi=(\varphi \theta _1^{-1})\cdot \text{diag}(y, 1, \cdots ,1)$. Thus the $x$-Wiebe matrix $\psi$ has the form described in (iv).
 
$(iv)\Rightarrow (ii)$ \quad  In the given factorization of the $x$-Wiebe matrix $\psi$, the first factor $\varphi _1$ is an $x$-nice matrix belonging to $(0:yA)=I_1$, see \ref{cwa} again, thus this ideal $I_1$ is  generated by the sequence $x\cdot \varphi _1$, say $x\cdot \varphi _1 = (z, a_2, \cdots ,a_n)$. But the sequence $x\cdot \psi=x\cdot \varphi _1 \cdot \text{diag}(y, 1, \cdots,1)= (yz, a_2,\cdots , a_n)$ is the zero sequence since $\psi$ is an $x$-Wiebe matrix. Thus we have $yz=0$ and  $a_i=0$ for $2\leq i\leq n$, in particular we have $I_1=zA$, $I_1$ is principal.

$(ii) \Rightarrow (iii)$ \quad If the ideal $(0:yA)$ is principal, its annihilator $yA$ is Gorenstein principal and thus C.I.0, see \ref{cku}.

$(iii) \Rightarrow (i)$ \quad This is contained in \ref{cku}.

$(iv) \Rightarrow (v)$ \quad In an $x$-realization $(S, Q, x')$ of $A$, see \ref{d1.2} for the notations, we take an element $y'$ which maps onto $y$, we also take a matrix $\varphi '_1 \in S^{n\times n}$ which maps onto the matrix $\varphi _1$ given by $(iv)$ and we  put $\psi '= \varphi '_1 \cdot \text{diag}(y', 1, \cdots , 1)$, so that the image of $\psi '$ in $A^{n\times n}$ is the $x$-Wiebe matrix $\psi$. With \ref{twi}(a), we know that the sequence $x'\cdot \psi '$ is a maximal regular sequence in $S$ generating the ideal $Q$. But $x'\cdot \psi '= x'\cdot \varphi '_1 \cdot \text{diag}(y', 1, \cdots , 1)$. Thus, if $x'\cdot \varphi '_1 = (z', a'_2, \cdots , a'_n)$, the ideal $Q$ is generated by the regular sequence $(y'z', a'_2, \cdots ,a'_n)$.

$(v) \Rightarrow (iii)$ \quad This a consequence of the following well-known fact (after a  permutation of our regular sequence).
\end{proof}

\begin{fact}\label{fr} In any commutative ring $A$ (with unit), if $(a_1, a_2, \cdots , a_{n-1}, yz)$ is a regular sequence generating an ideal $Q$, then both sequences $(a_1, \cdots ,a_{n-1}, y)$ and $(a_1, \cdots ,a_{n-1}, z)$ are regular, generating ideals $I_1$ and $I_2$ respectively, and we have $Q\subset I_1\cap I_2, \quad Q:I_1=I_2, \quad Q:I_2=I_1, \quad I_1/Q=(yA+Q)/Q, \quad I_2/Q=(zA+Q)/Q$. (Note that "regular" here can be replaced by "completely secant" and that this fact then may be viewed as a very particular case of \ref{p1.1}.)
\end{fact}

In the above theorem, the hypothesis that $y$ is a minimal generator of the maximal ideal cannot be dropped without harm, there are non principal C.I.0-ideals, but  their annihilator is in $M^2$.

\wl When we have in our ring $A$ a minimal generator of the maximal ideal $M$ which generates a C.I.0-ideal, it is a good idea to take it as the first element of a sequence $x=(x_1, \cdots , x_n)$ generating $M$.
Doing so, we have with (\ref{pgo},$(iii) \Leftrightarrow (iv)$) an $x$-Wiebe matrix for $A$ the first column of which is a multiple of $x_1$, but we also obtain another $x$-Wiebe matrix the first column of which is $\begin{pmatrix}z_1\\0\\ \vdots \\0 \end{pmatrix}$, where $(0:x_1A)=z_1A$. More precisely we have the following.

\begin{corollary}\label{cblabla} Let $A$ be a complete intersection of dimension zero with a sequence $x=(x_1, \cdots ,x_n)$ generating its maximal ideal $M$ minimally, and let us denote by $\overline{(\cdot )}$ the images modulo $x_1A$.

Then the following conditions are equivalent.

(i) The ideal $x_1A$ is C.I.0.

(ii) There is an $x$-Wiebe matrix for $A$ of the form 
$$\psi = \left( \begin{array}{c|c} z_1 & \cdots \\ \hline  0 & \psi^* \end{array}\right).$$

When these conditions are satisfied, let $\psi$ be a matrix as in (ii).\\
We then have:
 
 (a) $(0:x_1A)=z_1A$, 
 
 (b) $\overline{\psi^*}$ is an $(\overline {x}_2, \cdots , \overline {x}_n)$-Wiebe matrix for $\bar A=A/x_1A$.
\end{corollary}
\begin{proof} $(i) \Rightarrow (ii)$ \quad Let $(S,Q,x')$ be an $x$-realization of $A$. If $x_1A$ is a C.I.0-ideal, we know with \ref{pgo} that the ideal $Q$ of $S$ is generated by a maximal regular sequence of the form $(x'_1 z'_1, a'_2, \cdots , a'_n)$ which may be written
$$(x'_1z'_1, a'_2, \cdots , a'_n) = (x'_1, x'_2, \cdots , x'_n)\cdot \psi ', \text{ with } \psi '= \left( \begin{array}{c|c} z'_1 & \cdots \\ \hline  0 & \psi '^* \end{array}\right), \psi'^* \in S^{(n-1)\times (n-1)}.$$
We take the image $\psi$ of  $\psi '$ in $A$. With \ref{p1.4} we know that $\psi$ is an $x$-Wiebe matrix for $A$ and it has the wanted form.

 $(ii)\Rightarrow (i), (a), (b)$ \quad In the given $x$-Wiebe matrix $\psi$, we replace $z_1$ by 1, we obtain a matrix $\varphi _1$ for which $J(x\cdot \varphi _1) =x_1A$ and a factorization $\psi = \varphi _1 \cdot \text{diag}(z_1, 1, \cdots , 1)$. We conclude with \ref{cwa} that $x_1A$ is a C.I.0-ideal with $(0:x_1A)=z_1A$.
 
To finish, we take images modulo $x_1A$. We know with \ref{crw} that $\overline{\varphi} _1$  is an $(\overline{x}_1, \overline{x}_2, \cdots ,\overline{x}_n)$-Wiebe matrix  for $\bar A$ and we conclude that $\overline{\psi^*}$ is an  $(\overline{x}_2, \cdots ,\overline{x}_n)$-Wiebe matrix for $\bar A$.
\end{proof}

With \ref{pgo} we are also able to complete the information given in \ref{equi} about  strict chains of C.I.0-ideals of the maximal length and those C.I.0-ideals which can be taken as the starting point of such a chain.
 
\begin{corollary}\label{cnew}Let $A$ be a complete intersection of dimension zero and exponent $t+1$.

(i) \quad Let $I$ be a Gorenstein ideal with $\text{exponent}(A/I)= t$.Then:

$I$ is a C.I.0-ideal $\Leftrightarrow$ $I$ is principal $\Leftrightarrow$ $(0:I)$ is a principal C.I.0-ideal generated by some $y\in M\setminus M^2$.

(ii)\quad Assume we have in $A$ a strict saturated chain of C.I.0-ideals of the maximal length t: \quad $0\subset I_1\subset \cdots \subset I_t=M$. Then the following holds.

(a) $I_i/I_{i-1}$ is a principal C.I.0-ideal of $A/I_{i-1}$.

 The ideal $I_{i-1}:I_i$ is also a C.I.0-ideal of $A$ for all $i, 1\leq i \leq t$, more precisely $(I_{i-1}:I_i)/I_{i-1}$ is a principal C.I.0-ideal of the quotient ring $A/I_{i-1}$ generated by a minimal generator of its maximal ideal.

(b) In particular there is a sequence $z=(z_1, \cdots ,z_t)$ generating $M$ (not minimally if $\text{embedding dimension}(A)<t$) such that $I_i=z_1A+z_2A+\cdots z_iA$ , $1\leq i\leq t$.

And for any such sequence, there is an upper triangular $z$-Wiebe matrix for $A$ with main diagonal	 $\text{diag}(d_1, \cdots ,d_n)$ such that $I_{i-1}:I_i=I_{i-1}+d_iA, \quad  d_i \in M\setminus M^2$.

(iii) Conversely, assume we have a sequence $z=(z_1, \cdots ,z_t)$ generating the maximal ideal $M$ of $A$ and an upper triangular $z$-Wiebe matrix  with main diagonal $\text{diag}(d_1,  \cdots , d_t)$, where $d_i \in M, 1\leq i \leq t$.

Then the ideals $I_i=(z_1A+z_2A+ \cdots +z_iA)$ are C.I.0-ideals, they form a strict saturated chain and $I_{i-1}:I_i=I_{i-1}+d_iA$.
\end{corollary}

\begin{proof} (i) We know with (\ref{equi}(ii)(a)) that the ideal $0:I$ is principal, generated by a minimal generator of the maximal ideal $M$ of $A$, say $0:I=yA, y\in M\setminus M^2$. Thus $I=(0:yA)$ and we conclude with \ref{pgo}.

 (ii)(a) This follows from (i) applied to the ideal $I_i/I_{i-1}$  of $A/I_{i-1}$, since with \ref{equi} we know that $\text{exponent}(A/I_i) =\text{exponent}(A/I_{i-1})-1$.
 
 (ii)(b) The existence of the sequence $z$ such that $I_i=z_1A+ \cdots +z_iA$ is clear from (ii)(a).
 
 We take images modulo $I_1=z_1A$ and denote them by $\overline{(\cdot)}$, so $\bar A=A/I_1$. An induction on the exponent gives us an upper triangular matrix $\psi^*\in A^{(t-1)\times (t-1)}$ such that $\overline{\psi^*}$ is a $(\overline{z}_2, \cdots , \overline{z}_t)$-Wiebe matrix for $\bar A$ of the wanted form, which means that $J((z_2, \cdots , z_t)\cdot \psi^*)\subset z_1A$, that $\text{det}(\psi^*)\notin z_1A$ and that, if $\text{diag}(d_2, \cdots , d_t)$ is the main diagonal of $\psi^*$, we have $I_{i-1}:I_i = I_{i-1}+d_iA, \quad 2\leq i\leq t$.
 
 Now let $\beta =\left( \begin{array}{c|c} d_1&0\\ \hline 0&\psi^*\\ \end{array} \right)$, where $(0:I_1)=d_1A$. It is easy to see that $J(z\cdot \beta)\subset z_1A$. This implies that we have a sequence $(a_2, \cdots ,a_n)$ such that the matrix $$\psi= \left( \begin{array}{c|c} d_1& a_2 \cdots a_n\\
 \hline 0&\psi^*
 \end{array}\right)$$
 satisfies $J(z\cdot \psi)=0$. As $\text{det}(\psi^*) \notin z_1A=0:d_1A$, we also have $\text{det}(\psi)=d_1 \text{det}(\psi ^*) \neq 0$, thus $\psi$ is the $z$-Wiebe matrix we were looking for.
 
 (iii) Let $\psi$ be an upper triangular $z$-Wiebe matrix with main diagonal $\text{diag}(d_1, \cdots ,d_n)$, where the $d_i$'s are in $M$. Since $\text{det}(\psi)= d_1\cdots d_t$ generates the socle $M^t$ of $A$, we first note that $d_i \in M\setminus M^2$ for all $i, 1\leq i \leq t$.
 
 We factorize $\psi$:
 $$\psi = \varphi _1 \cdot \text{diag}(d_1, 1, \cdots ,1), \text{ where } \varphi _1=\left( \begin{array}{c|c} 1&\cdots \\ \hline 0& \psi ^* \end{array}\right).$$
 We observe that $J(z\cdot \varphi _1)=z_1A$, we write $I_1=z_1A$ and we conclude with \ref{cwa} that $I_1$ is a C.I.0-ideal, that $0:I_1=d_1A$.
 
 We now take images modulo $z_1A$ and denote them by $\overline {(\cdot )}$.
 
 From $\text{det}(\psi) =d_1 \text{det}(\psi ^*)\neq 0$, we deduce that $\text{det}(\psi ^*) \notin z_1A$ and we then observe that $\overline{\psi^*}$ is an upper triangular  $(\overline z _2, \cdots , \overline z _t)$-Wiebe matrix for $\bar A$.
 
 Since $\text{exponent}(\bar A)< \text{exponent}(A)=t+1$ and $\text{det}(\overline{\psi ^*}) \in M^{t-1}$, we also have $\text{exponent}( \bar A)=t$.
 
 We now conclude with an induction on the exponent.
\end{proof}

Later on, we shall also see  in \ref{c5t} some  strange conditions under which a principal ideal generated by a minimal generator of the maximal ideal is C.I.0.

Here are some other sufficient conditions.

\begin{lemma}\label{ltri}Let $A$ be a complete intersection of dimension zero with maximal ideal $M$ and assume we have in $A$ two elements $y,z \in M\setminus M^2$ such that $yz=0$.

Then the principal ideals $yA$ and $zA$ are both C.I.0-ideals and
$0:yA=zA,\quad 0:zA=yA$.

If moreover $A$ is positively graded and if  $y$, $z$ are homogeneous, then $\text{exponent}(A/yA)=\text{exponent}(A/zA)=\text{exponent}(A)-1$.
\end{lemma}
\begin{proof}  Let $x=(x_1, \cdots , x_n)$ be a sequence generating the maximal ideal of $A$ minimally, let $(S, Q, x')$ be an $x$-realization of $A$ (see \ref{d1.2} for the notations), and let $y', z' \in S$ be such that their images in $A$ are $y, z$ respectively. We denote as usual by  $M'$ the maximal ideal of $S$. 

We have $y', z' \in M'\setminus M'^2$ and $y'z' \in Q$ by construction, $Q\subset M'^2$ by the minimality of the sequence $x$, and we also have $y'z' \in M'^2\setminus M'^3$ since $M'$ is generated by a regular sequence. Thus $y'z'$ is a minimal generator of $Q$. On the other hand we know with \ref{twi} that $Q$ is generated by a regular sequence,  we then know that every sequence generating $Q$ minimally is regular. Thus $Q$ is generated by a regular sequence of the form $(a'_1, a'_2, \cdots ,a'_{n-1}, y'z')$. We conclude with \ref{fr}.

The second assertion follows from the first and (\ref{equi} (ii)(b)).
\end{proof}

We now turn to a particular class of complete intersections (where we have some hope to encounter the conditions of \ref{ltri}). The following terminology is justified by \ref{c1.3}.

\begin{definition}\label{dminex}
We shall say that a complete intersection $A$ of dimension zero has {\bf minimal exponent} if $\text{exponent}(A)=\text{embedding dimension}(A)+1$. 
\end{definition}

For the rings in this particular class, the sufficient condition in \ref{ltri} is also necessary.

\begin{proposition}\label{pimi} Let $A$ be a complete intersection of dimension zero and minimal exponent $n+1$ and let $y$ be a minimal generator of its maximal ideal $M$.

Then $yA$ is a C.I.0-ideal if and only if there is an element $z\in M\setminus M^2$such that $yz=0$.

In that case $(0:yA)=zA$.
\end{proposition}
\begin{proof} The if part and the last assertion is contained in \ref{ltri}.

Assume now that $yA$ is a C.I.0-ideal. Since the nonnull ideal $(0:yA)$ contains the socle $(0:M)=M^n$ of $A$, we have $\text{exponent}(A/(0:yA))\leq n$. Since $(0:yA)$ is also a C.I.0-ideal, see \ref{pgo}, we then have with \ref{c1.3} $\text{embedding dimension}(A/(0:yA))< n$, which means that $(0:yA)$ contains a minimal generator of $M$.
\end{proof}

We now complete the information given in \ref{pgo} and \ref{cnew}.

\begin{proposition}\label{cmin} Let $A$ be a complete intersection of dimension zero and minimal exponent $n+1$.
 
 (i)\quad Assume we have in $A$ a C.I.0-ideal $I$ such that $\text{exponent}(A/I)=\text{exponent}(A)-1$.
 
 Then both ideals $I$ and $(0:I)$ are principal C.I.0-ideals and each of them is generated by a minimal generator of the maximal ideal $M$.
 
(ii) There is in $A$ a strict saturated chain of C.I.0-ideals of the maximal length $n$ if and only if there is a sequence $z=(z_1, \cdots , z_n)$ generating the maximal ideal $M$ of $A$ minimally and an upper triangular $z$-Wiebe matrix for $A$.. 
 \end{proposition}

\begin{proof} 

(i) \quad With \ref{cnew} we already know that $I$ and $(0:I)$ are principal C.I.0-ideals and that $(0:I)=yA$ for some $y\in M\setminus M^2$. With \ref{pimi} we then have that $I$ contains a minimal generator $z$ of $M$ and we conclude that $I=zA$.
 
(ii) \quad This is a direct consequence of \ref{cnew}, remember that, when a sequence $z$ generates $M$ minimally, all the entries of a $z$-Wiebe matrix belong to $M$.
\end{proof}

\begin{observation} \label{oz} If $A$ is a complete intersections of dimension zero and  minimal exponent $n$, if $x=(x_1, \cdots , x_n)$ is a sequence generating the maximal ideal $M$ of $A$ minimally and if $(S, Q, x')$ is an $x$-realization of $A$, we  note that the ideal $Q$ of $S$ is generated by elements in $M'^2\setminus M'^3$, where $M'$ denotes the maximal ideal of $A$.

Indeed, $Q$ is generated by a maximal regular sequence which may be written $(a'_1, \cdots , a'_n) = x'\cdot \psi '$ for some matrix $\psi ' \in S^{n\times n}$. Since $Q \subset M'^2$ by the minimality of the sequence $x$, we see that all the entries of any such  $\psi '$ belong to $M'$. If one of these generators $a'_i$ of $Q$ is in $M'^3$, then we can take $\psi '$ such that  one of its columns  has its entries in $M'^2$. But the image of this  $\psi '$ in $A$ is an $x$-Wiebe matrix $\psi$ for $A$ with 
$\text{det}(\psi)A=0:M=M^n\neq 0$. This implies that each column of $\psi$ must have at least one entry outside $M^2$.
\end{observation}

\section{Towards a classification of C.I.0-ideals}

In this  section $A$ is again a noetherian local ring with a sequence $x=(x_1,\cdots ,x_n)$ generating its maximal ideal $M$, and we assume that $A$ is complete in its $M$-adic topology, so that it has an $x$-realization.

 The set of $n\times n$ matrices with entries in $A$, with its ring structure, will be denoted by ${\cal M}_n(A)$. Again we often identify a map $A^m \rightarrow A^n$ with the $n\times m$ matrix representing it on the canonical bases.

Now let $\varphi _1$, $\varphi _2 \in {\cal M}_n(A)$. When do we have the inclusion $J(x\cdot \varphi _1) \subset J(x\cdot \varphi _2)$, when do we have the equality?

Of course, we have the inclusion if and only if there is a matrix  $\gamma$ in ${\cal M}_n(A)$ such that
 $x\cdot ( \varphi _1 - \varphi _2 \cdot  \gamma) = 0$, but we can say more,
 at least when one of the matrices involved  is $x$-nice. Before, we need some preliminaries.

\begin{definition}\label{d41} Let $\delta$ be the matrix representing the second boundary map in the Koszul complex $K_{\cdot}(x, A)$, so that  $ \delta \in A^{n\times \frac{n(n-1)}{2}}, \quad  x \cdot \delta  =0$. By the {\bf $x$-Koszul ideal} we mean the right-ideal in the ring ${\cal M}_n(A)$ defined by
 $${\cal N}_{x,A} =  \{ \alpha \in {\cal M}_n(A) \mid \alpha \phantom{i}\text{can be factored through } \delta\} =\{ \alpha \in {\cal M}_n(A) \mid \text{im}(\alpha) \subset \text{im}(\delta) \}.$$

\[
\xymatrix{
A{\frac{n(n-1)}{2}} \ar[rr]^-{\delta} \hole&&  A^n  \ar[rr]^- {x} \hole&& A\\
 & A^n \ar@{.>}[ul]^{\beta} \ar[ur]_ {\alpha} }
 \]
 Here of course we have replaced matrices by the maps they induce on the appropriate free bases wherever convenient.
\end{definition}

In general, the $x$-Koszul ideal is smaller than the right-annihilator ideal of $x$ in ${\cal M}_n(A)$, these two ideals coincide exactly when the Koszul complex $K_{\cdot}(x,A)$ is exact in degree one, i.e. when the sequence $x$ generating the maximal ideal of $A$ is regular (which means that the sequence $x$ generate minimally the maximal ideal of $A$ and  that the ring $A$ is regular).However, this $x$-Koszul ideal is enough for our purpose, as we shall see.

\begin{remark}\label{rb}We note that an $n\times n$ matrix belongs to the $x$-Koszul ideal if and only if each column of it belongs to the image of the second boundary map $\delta$ in the Koszul complex.

With a little reflection, we also see that, if $\bar A$ is a homomorphic image of $A$ and $\bar x$ is the image of $x$ in $\bar A$, then ${\cal N}_{x,A}$ maps onto   ${\cal N}_{\bar x,\bar A}$~. (Indeed, if $\overline \alpha =  \overline \delta \cdot \overline \beta$ belongs to ${\cal N}_{\bar x,\bar A}$, then $\delta \cdot \beta$ belongs to ${\cal N}_{x, A}$ and maps onto $\overline \alpha$.)
\end{remark}
 
 \begin{proposition}\label{p4in}Let $\varphi _1 , \varphi _2 \in {\cal M}_n(A)$ and assume that the matrix $\varphi _2$ is $x$-nice.
 
 Then \quad $J(x\cdot \varphi _1) \subset J(x\cdot \varphi _2)$ \quad if and only if \quad $\varphi _1 \in  \varphi_2 \cdot {\cal M}_n(A) + {\cal N}_{x,A}$. 
  \end{proposition}
  
 \begin{proof}
 The if part is obvious. So assume that  $J(x\cdot \varphi _1) \subset J(x\cdot \varphi _2)$. We take an $x$-realization $(S,Q,x')$ of $A$, $A \simeq  S/Q$, and we take matrices $\varphi '_1, \varphi' _2 $ in ${\cal M}_n(S)$ the images of which in ${\cal M}_n(A)$ are $\varphi _1$ and $\varphi _2$ respectively. Thus $J(x'\cdot \varphi' _1) \subset J(x' \cdot \varphi '_2) + Q$. 
 But $Q \subset J(x' \cdot \varphi ' _2)$ since the matrix $\varphi _2$ is $x$-nice, see lemma \ref{l}, we  thus have $J(x' \cdot \varphi ' _1)\subset J(x' \cdot \varphi '_2)$ and there is a matrix $\gamma' \in {\cal M}_n(S)$ such that $x'\cdot \varphi ' _1 = x' \cdot \varphi ' _2 \cdot \gamma'$. This implies that $\varphi ' _1 - \varphi ' _2 \cdot \gamma ' \in {\cal N}_{x',S}$ since the sequence  $x'$ is regular on $S$.
 As the image of ${\cal N}_{x',S}$ in ${\cal M}_n(A)$ is ${\cal N}_{x,A}$, we obtain what we wanted.
 \end{proof} 
  
   One point is in order and the preceding leads us to a new definition
 
 \begin{definition}\label{dns} A submodule of the right ${\cal M}_n(A)$-module ${\cal M}_n(A)/ {\cal N}_{x,A}$ is said to be {\bf $x$-nice} if it is  cyclic  and can be generated by the image of an $x$-nice matrix.
  \end{definition}
  
  With the above and the characterization of C.I.0-ideals given in corollary \ref{crw}, we have the following.
  
 \begin{theorem}\label{cfirst} (i) Let  $\varphi _1 , \varphi _2 \in {\cal M}_n(A)$ be two $x$-nice matrices.
 
 Then these matrices are equivalent, i.e.  $J(x\cdot \varphi _1) = J(x\cdot \varphi _2)$
 
  \quad $\Leftrightarrow \quad \varphi _1 \cdot {\cal M}_n(A) + {\cal N}_{x,A} =  \varphi_2 \cdot {\cal M}_n(A) + {\cal N}_{x,A}$.
 
  (ii) The map from the set $\{ \varphi \in {\cal M}_n(A) \mid \varphi \text{ is x-nice} \}$ to the set of  $x$-nice submodules of ${\cal M}_n(A)/{\cal N}_{x,A}$ given by $$\varphi_1 \mapsto  (\varphi _1 \cdot {\cal M}_n(A) + {\cal N}_{x,A})/{\cal N}_{x,A}$$ induces a bijection preserving inclusion between the set of C.I.0-ideals of $A$ and the set of  $x$-nice submodules of ${\cal M}_n(A)/{\cal N}_{x,A}$.
  \end{theorem} 
  
  Our next aim is the study of  all generators of a given $x$-nice submodule of ${\cal M}_n(A)/{\cal N}_{x,A}$, are they all image of an $x$-nice matrix?  
 
  And, in analogy with \ref{pfac}, is it true  that a submodule of ${\cal M}_n(A)/{\cal N}_{x,A}$ containing an $x$-nice one is  also $x$-nice? To provide a positive answer we need again some preliminaries.
  
 \begin{lemma}\label{ll4} Let $\varphi$ be an $x$-nice matrix. Then
 
 (i) \quad $ \forall \alpha \in {\cal N}_{x,A}$, \quad $\varphi + \alpha$ is an $x$-nice matrix (equivalent to $\varphi$).
 
 (ii) If $\gamma \in {\cal M}_n(A)$ is such that $\varphi - \varphi \cdot \gamma \in {\cal N}_{x,A}$, then $\gamma$ is invertible.
 \end{lemma}
  
 \begin{proof} Let again $(S,Q,x')$ be an $x$-realization of $A$, $A \simeq  S/Q$, and let $\varphi ', \gamma '$ in ${\cal M}_n(S)$  be two matrices, the images of which in ${\cal M}_n(A)$ are $\varphi$ and $\gamma$ respectively, let also $\alpha ' \in {\cal N}_{x',S}$ be a matrix the image of which is $\alpha$ (the existence of such an $\alpha '$ follows from \ref{rb}).
 
(i) \quad As $\varphi$ is $x$-nice, we have that $Q\subset J(x'\cdot \varphi ')$ and that the matrix $\varphi '$ is $x'$-nice, see  lemma \ref{l}, so that the sequence $x'\cdot \varphi '$ is a maximal regular sequence in $S$, see proposition \ref{p1.4}.  But $x'\cdot \varphi ' = x' \cdot (\varphi ' + \alpha ')$, so \ref{p1.4} again gives us $\text{det}(\varphi ' + \alpha ') \notin J(x' \cdot (\varphi ' + \alpha')) = J(x'\cdot \varphi ')$. Going back into $A$ we obtain $\text{det}(\varphi + \alpha ) \notin J(x \cdot (\varphi + \alpha )) $.

(ii) \quad Assume now that $\varphi - \varphi \cdot \gamma \in {\cal N}_{x,A}$. Then $\varphi \cdot \gamma$ is $x$-nice equivalent to $\varphi$ by (i) and $J(x\cdot \varphi) = J(x\cdot \varphi\cdot \gamma)$, so that $Q + J(x'\cdot \varphi ') = Q + J(x'\cdot \varphi ' \cdot \gamma')$.
 As above, since the matrices $\varphi $ and $\varphi \cdot \gamma $ are $x$-nice,  we also have  $J(x'\cdot \varphi ') = J(x'\cdot( \varphi '\cdot \gamma '))$ and  both sequences  $x'\cdot \varphi '$ and $x'\cdot \varphi '\cdot \gamma '$ are maximal regular sequences in $S$  generating the ideal $J(x'\cdot \varphi ')$ minimally. But any matrix obtained by writing the elements of a minimal set of generators of a finitely generated $A$-module as linear combinations of  the elements of another minimal set of generators of this module is invertible because its determinant is invertible. Thus $\gamma '$ is invertible and so is $\gamma$.
 \end{proof}
 
 \begin{corollary}\label{ct} Let $\beta \in {\cal M}_n(A)$. If any single column of $\beta$ belongs to the image of the second boundary map $\delta$ in the Koszul complex $K_{\cdot}(x, A)$, then the matrix $\beta$ is not $x$-nice.
 \end{corollary}
 \begin{proof} We can write $\beta = \beta _1 + \alpha$, where $\alpha \in {\cal N}_{x,A}$ and $\beta _1$ is a matrix, one column of which is null. As $\text{det}(\beta _1) =0$, $\beta _1$ is not $x$-nice and so is $\beta$.
 \end{proof}
 
 The converse of the above corollary is false. Here is an example.
 
 \begin{example}\label{r4}
 
 Let $A= K[[X,Y]]/(X^2 + Y^2, XY)$, where $K$ is a field, and let us denote as usual by $x,y$ the images of $X,Y$ in $A$.
 This ring $A$ is a complete intersection of dimension zero with $(x,y)$-Wiebe matrix $\psi = \left( \begin{array}{cc} x&0\\y&x \end{array} \right)$.
 
 Now let $\beta =\left( \begin{array}{cc}
 x&0\\
 0&y
 \end{array}\right)$, we have $\text{det}(\beta)=xy=0$, the matrix $\beta$ is not $(x,y)$-nice though its columns do not belong to $\text{im}(\delta) =  \left( \begin{array}{c} -y\\ x \end{array} \right) A $.
 \end{example}
 
 Here are some answers to the questions raised before \ref{ll4}.
 
 \begin{proposition}\label{pin} Let $\varphi _1, \varphi _2 \in {\cal M}_n(A)$. 
 
 If the matrix $\varphi _1$ is $x$-nice and if $\varphi _1 \in \varphi _2 \cdot {\cal M}_n(A) + {\cal N}_{x,A}$, then the matrix $\varphi _2$ is also $x$-nice and and $J(x\cdot \varphi _1) \subset J(x\cdot \varphi _2)$.
 \end{proposition}
 
 \begin{proof} We write $\varphi _1 = \varphi _2 \cdot \gamma + \alpha$, where $\gamma \in {\cal M}_n(A)$ and $\alpha \in {\cal N}_{x,A}$.
 With \ref{ll4} we know that the matrix $ \varphi _2 \cdot \gamma$ is $x$-nice and we obtain that $\varphi _2$ is also -x-nice with proposition \ref{pfac}. The inclusion is obvious.
 \end{proof}

 \begin{corollary}\label{c42} Let $I$ be a C.I.0-ideal of the ring $A$ and let $\varphi$ be an $x$-nice matrix belonging to it: $I=J(x\cdot \varphi)$. Then
 
 \wl (i) the equivalence class of $\varphi$, i.e. the set of all $x$-nice matrices belonging to $I$ is the set
  $$\{\varphi\cdot \theta + \alpha \mid \alpha \in {\cal N}_{x,A},\quad \theta \in {\cal M}_n(A), \quad \theta \quad \mbox{is invertible}\}$$
  
  (ii) this equivalence class is also equal to $$\{ \varphi ' \in {\cal M}_n(A) \mid \varphi \cdot {\cal M}_n(A)+{\cal N}_{x,A} = \varphi ' \cdot {\cal M}_n(A)+{\cal N}_{x,A} \}$$ 
  In particular, every generator of an $x$-nice submodule of the right ${\cal M}_n(A)$-module ${\cal M}_n(A)/{\cal N}_{x,A}$ is the image of an $x$-nice matrix.
 \end{corollary}
 
 \begin{proof}  (i) With the lemma we have that the matrices of the form $\varphi \cdot \theta + \alpha$, where $\theta \in {\cal M}_n(A)$ is invertible and $\alpha \in {\cal N}_{x,A}$, are $x$-nice matrices belonging to $I$.
 
 Conversely, let $\varphi _1$ be another $x$-nice matrix belonging to $I$. With theorem \ref{cfirst}, we obtain $\varphi \cdot {\cal M}_n(A) + {\cal N}_{x,A} =  \varphi_1 \cdot {\cal M}_n(A) + {\cal N}_{x,A}$ and we can write 
 
 $\varphi _1 =\varphi \cdot \theta + \alpha,\quad \varphi = \varphi _1 \cdot \theta _1 + \alpha _1$, where $\theta, \theta_1 \in {\cal M}_n(A)$ and $\alpha, \alpha _1 \in {\cal N}_{x,A}$.\\
 We then have $\varphi _1 = (\varphi _1 \cdot \theta _1 + \alpha _1)\cdot \theta + \alpha$ and $\varphi _1 - \varphi _1 \cdot \theta _1 \cdot \theta \in {\cal N}_{x,A}$.
 Now  lemma \ref{ll4} tells us that the matrix $\theta _1 \cdot \theta$ is invertible and so is $\theta$.
 
  (ii) A matrix $\varphi '$ as in (ii) is necessarily $x$-nice by the previous proposition and is an $x$-nice matrix belonging to $I$, this gives us one inclusion. The other inclusion is theorem \ref{cfirst}. 
 \end{proof}

 The above proposition \ref{pin} also give us the following.

 \begin{corollary}\label{t4} Any cyclic right submodule of ${\cal M}_n(A)/{\cal N}_{x,A}$ containing an $x$-nice submodule is $x$-nice.
 \end{corollary}

 The case when the ring $A$ itself is a complete intersection of dimension zero deserves attention.
 
 \begin{theorem}\label{clast} Assume that the ring $A$ is a complete intersection of dimension zero and let ${\cal M}_0$ be the $x$-nice submodule of ${\cal M}_n(A)/{\cal N}_{x,A}$ corresponding to the null ideal of $A$, thus ${\cal M}_0 = (\psi\cdot {\cal M}_n(A) + {\cal N}_{x,A})/{\cal N}_{x,A}$, where $\psi$ is an $x$-Wiebe matrix for $A$.
 
 Then the C.I.0-ideals of $A$ are in bijection with the cyclic submodules of \quad ${\cal M}_n(A)/ {\cal N}_{x,A}$ containing ${\cal M}_0$.
 \end{theorem}

\section{Last remarks}
In a regular local ring $S$ with a sequence $x'=(x_1', \cdots ,x_n')$ generating its maximal ideal minimally, the $x'$-nice matrices correspond to the maximal regular sequences in $S$ and give us a description of the elements of the corresponding sequence, see \ref{p1.4}.  In particular any $x'$-nice matrix $\varphi '$ belonging to a C.I.0-ideal $I'$ of $S$ gives us a minimal set of generators of the ideal $I'=J(x'\cdot \varphi ')$, the matrix $\varphi '$ encodes the structure of the zero-dimensional complete intersection $S/J(x'\cdot \varphi ')$.

In the same way, in a ring $A$ with a sequence $x=(x_1, \cdots ,x_n)$ generating its maximal ideal, an $x$-nice matrix encodes the structure of the quotient $A/J(x\cdot \varphi)$, see \ref{rwi}.

In view of this the following has perhaps some interest. At least, it will provide conditions under which the ideal generated by a minimal generator of the maximal ideal is a C.I.0-ideal.

\begin{remark}\label{r51}Let $x$ be a sequence generating the maximal ideal of a complete noetherian local ring $A$. In view of corollary \ref{c42}, we may perform some column operations on an $x$-nice matrix $\varphi$ without affecting its equivalence class, which means that the matrix $\varphi _1$ obtained after these operations is sill an $x$-nice matrix belonging to $J(x\cdot \varphi)$. 

First we have what we call the column  operations, we may permute two columns of an  $x$-nice matrix $\varphi$, we may multiply one column of it by an invertible element of the ring, we may add to one of its  columns a scalar multiple of another column of the matrix. These operations amount to multiply our matrix on the right by an invertible matrix and do not affect its equivalence class.

We may also play with the image of the second boundary map $\delta$ in the Koszul complex $K_{\cdot}(x, A)$. We may add a column matrix belonging to $\text{im}(\delta)$ to one of the columns of an $x$-nice matrix $\varphi$, doing so  we still have an $x$-nice matrix equivalent to the one we started with, see \ref{ll4}. Indeed, this amount to replace the matrix $\varphi$ by a matrix of the form $\varphi + \alpha$, where $\alpha \in {\cal{N}}_{x,A}$.
\end{remark}

 \begin{lemma}\label{l52} Let $x=(x_1,\cdots ,x_n)$ be a sequence generating the maximal ideal of a complete noetherian local ring $A$.  An $x$-nice matrix $\varphi$ is always equivalent to an  $x$-nice matrix  of the form 
$$\varphi _1 = \left( \begin {array}{c|c}
x_1^{r_1} & 0\cdots 0 \\
\hline
\vdots & \phi^* _1 
\end{array} \right), \quad r_1\geq 0.$$

Moreover, let us put $H_{1,x}= (x_2A+ \cdots +x_nA )$. For such a matrix $\varphi_1$, the following holds. 

(i)\quad The natural number $r_1$  is determined by $$r_1= \text{max}\{i\in \mathbb N \mid x_1^i \notin ( H_{1,x}+ J(x\cdot \varphi))\}$$.

(ii)\quad   $(J(x\cdot \varphi):\text{det}(\phi^*_1)A)= (H_{1,x}+J(x\cdot \varphi))$.
 \end{lemma}
 \begin{proof}We note that each entry of the matrix $\varphi$ is a polynomial expression in the $x_i$'s with invertible coefficients. We recall that, $\forall i\geq 2$, $\text{im}(\delta)$ contains a column matrix whose first entry is $(-x_i)$, whose $i^{th}$ entry is $x_1$, the other entries being null. Thus, for all $y\in H_{1,x}= (x_2A+\cdots x_nA)$, there is in $\text{im}(\delta)$ a column matrix with $y$ as its first entry.
 
With \ref{r51} and the above, we  may add to each column of $\varphi$ a column belonging to $im(\delta)$ in order to obtain a matrix equivalent to $\varphi$ of the form
  $$\left (\begin{matrix}
u_1 x_1^{s_1} & u_2 x_1^{s_2} & \cdots & u_n  x_1^{s_n}\\
\hdotsfor[2]{4}\\
\hdotsfor[2]{4}\\
\hdotsfor[2]{4}\\
\end{matrix} \right ),\quad \begin{array}{c}\text{ where} \phantom{i}  u_1, u_2, \cdots, u_n \\ \text{are invertible or null}\\ \text{and where} \quad s_1,s_2,\cdots ,s_n\geq 0 .\end{array}$$

To obtain the wanted form $\varphi _1$, we just have to perform on this last matrix some column operations which amount to multiply it on the right by an invertible matrix.

\wl By definition the sequence $x\cdot \varphi _1$ generates the ideal $J(x\cdot \varphi _1) =J(x\cdot \varphi)$, so the relation $x_1^{r_1+1} \in   (H_{1,x}+ J(x\cdot \varphi ))$ is clear. Assume that we also have  $x_1^{r_1} \in  (H_{1,x} +J(x\cdot \varphi))$. Then we can add to the first column of $\varphi_1$ a column matrix belonging to $\text{im}(\delta)$ in order to obtain another $x$-nice matrix $\beta$ belonging to $J(x\cdot \varphi)$ and such that its (1,1) entry satisfies $\beta_{1,1}\in J(x\cdot \varphi)$. But then $\text{det}(\beta) = \beta _{1,1} \cdot \text{det}(\phi_1^*) \in  J(x\cdot \varphi)= J(x\cdot \beta)$, which is excluded for an $x$-nice matrix. This proves (i).

\wl We now factorize the matrix $\varphi _1$
$$\varphi _1=\varphi _2 \cdot \gamma, \quad \text{where}\quad \varphi _2 =\text{diag}(x_1^{r_1}, 1, \cdots ,1), \quad \gamma=\left(\begin{array}{c|c} 1&0\cdots 0\\ \hline  \vdots & \phi_1^* \end{array}\right).$$ 
We observe that $J(x\cdot \varphi _2)=(H_{1,x}+J(x\cdot \varphi))$, we recall that $J(x\cdot \varphi)= J(x\cdot \varphi _1)$ and we conclude with \ref{pfac} that $(J(x\cdot \varphi):\text{det}(\phi^*_1)A)= (H_{1,x}+J(x\cdot \varphi))$.
\end{proof}

  \wl We may apply the lemma to an $x$-Wiebe matrix. Here is a particular case.
 
 \begin{corollary}\label{5tr}Assume the ring $A$ is a complete intersection of dimension zero, embedding dimension $n\geq 2$ and minimal exponent $n+1$, which means that its maximal ideal satisfies  $M^n\neq 0,\quad  M^{n+1}=0$, and let $x=(x_1,\cdots , x_n)$ be any sequence generating the maximal ideal $M$ minimally. Then
 
 (i) \quad $A$ has an $x$-Wiebe matrix of the form $$\psi _1  =\left( \begin {array}{c|c}
x_1 & 0\cdots 0 \\
\hline
\vdots & \psi^* _1 
\end{array} \right) .$$

(ii)\quad  $x_1^2 \in (x_2A+ \cdots x_nA)$.

(iii)\quad  $(0:\text{det}(\psi _1^*)A) = (x_2A+\cdots +x_nA)$.

\end{corollary}
 
 \begin{proof} The entries of any $x$-Wiebe matrix $\psi$ for $A$  belong to $M$ since the sequence $x$ generates $M$ minimally, and the determinant of an $x$-Wiebe matrix generates the socle $M^n=(0:M)$ of $A$.  Thus, for the $x$-Wiebe matrix $\psi _1$ given by \ref{l52} we must have $\text{det}(\psi^*_1)\in M^{n-1}$, so that  $r_1=1$. 
 
 The other assertions follow now from \ref{l52}.
 \end{proof}

 When we looked at chains of C.I.0-ideals in a complete intersection of dimension zero, we already came across questions concerning the minimal generators of the maximal ideal, see \ref{pgo}, \ref{q4} . 
 
 Now the above lemma will give us sufficient conditions on a minimal generator of the maximal ideal to generate a C.I.0-ideal.

 \begin{proposition}\label{c5t} Let $A$ be a complete intersection of dimension zero and let $x_1$ be a minimal generator of the maximal ideal $M$ of $A$.
 
  If, for some sequence $x=(x_1,\cdots ,x_n)$ starting at $x_1$ and generating $M$ minimally, we have $$(0:x_1A) \nsubseteq (x_2A+ \cdots +x_nA) ,$$ then the principal ideal $x_1A$ is a C.I.0-ideal.
 \end{proposition}
 
 \begin{proof} Assume that the ideal $x_1A$ is not C.I.0 and let $x=(x_1,\cdots ,x_n)$ be any sequence starting at $x_1$ and generating $M$ minimally.
 
 With  lemma \ref{l52}, we have an $x$-Wiebe matrix for $A$ of the form 
 $$\psi =\left( \begin {array}{c|c}
x_1^{r_1} & 0\cdots 0 \\
\hline
\vdots & \psi^* 
\end{array} \right) .$$
We replace the first column of $\psi$ by $\left(\begin{array}{c} 1\\0\\ \vdots\\0\end{array}\right)$, we obtain a new matrix $\beta$ for which $J(x\cdot \beta)=x_1A$ (remember that $\psi$ was an $x$-Wiebe matrix, $J(x\cdot \psi)=0$) . With our assumption on $x_1A$, this matrix $\beta$ is not $x$-nice, thus $\text{det}(\beta)=\text{det}(\psi^*) \in x_1A$. We take annihilators and conclude with \ref{l52}: $(0:x_1A)\subseteq (0:\text{det}(\psi^*)A) = (x_2A+\cdots x_nA)$. 
 \end{proof}
 
 Here is a partial converse of the above proposition.
 
 \begin{corollary}\label{bof} Let $A$ be a complete intersection of dimension zero and minimal exponent $n$, which means that its maximal ideal satisfies $M^n\neq 0, M^{n+1}=0$, and let $x_1$ be a minimal generator of the maximal ideal.
 
 Then the ideal $x_1A$ is a C.I.0-ideal if and only if there is a sequence $x=(x_1,\cdots ,x_n)$ starting at $x_1$ and generating $M$ minimally such that $(0:x_1A) \nsubseteq (x_2A+ \cdots +x_nA) $.
 \end{corollary}
\begin{proof} The if part is \ref{c5t}.

Assume now that the ideal $x_1A$ is C.I.0. With \ref{pimi} we know that $0:x_1A= zA$ for some $z\in M\setminus M^2$, so that the images of $x_1$ and $z$ in the vector space $M/M^2$ over $A/M$ are nonnull. We distinguish two cases.

If $z\in x_1A +M^2$, then, for any sequence $x=(x_1, \cdots ,x_n)$ starting at $x_1$ and generating $M$ minimally, we have $z\notin (x_2A+\cdots x_nA)$.

If $z\notin x_1A +M^2$, then $(x_1, x_1+z)$ is a subsequence of a sequence $x=(x_1, x_1+z, x_3, \cdots ,x_n)$ generating $M$ minimally and for such a sequence $x$ we have $z\notin ((x_1+z)A+x_3A+\cdots x_nA)$.
\end{proof}

\

\

\begin{tabular}{ll}
Anne-Marie Simon                        &       Jan R. Strooker\\
Service d'Algebre C.P. 211              &       Mathematisch Instituut\\ 
Universite Libre de Bruxelles\qquad\qquad\qquad &       Universiteit Utrecht\\
Campus Plaine                           &       Postbus 80010\\
Boulevard du Triomphe                   &       3508 TA Utrecht,\\ 
B-1050 Bruxelles, Belgique              &       The Netherlands\\
e.mail: amsimon@ulb.ac.be               &       e.mail: strooker@math.uu.nl     
\end{tabular}

\end{document}